\documentclass[a4paper,11pt]{article}
\usepackage{amssymb,amsmath,latexsym,epsf,graphicx}

\newtheorem{thm}{Theorem}[section]
\newtheorem{lem}[thm]{Lemma}

\newtheorem{cor}[thm]{Corollary}
\newtheorem{prop}[thm]{Proposition}

\newtheorem{rem}[thm]{Remark}

\title{Counting invertible Schr\"odinger Operators over Finite Fields for Trees, Cycles and Complete Graphs.}
\author{Roland Bacher\footnote{This work has been partially supported by the LabEx PERSYVAL-Lab (ANR--11-LABX-0025). The author is a member of the project-team GALOIS supported by this LabEx.}}

\begin{document}
\maketitle

{\sl Abstract\footnote{
Math. class: Primary: 05C31, Secondary: 05C10, 20G40}: 
We count invertible Schr\"odinger operators (perturbations by diagonal 
matrices of the adjacency matrix) over finite fields
for trees, cycles and complete graphs.}

%fichier schroedtree1.tex dans recherche/schroedtree
%http://edf.lms.ac.uk/ef/status.php?p_id=26399&cr=748C1E1B85
\section{Introduction}

A \emph{Schr\"odinger
operator} (or perhaps more accurately, 
an opposite of a Schr\"odinger operator) on a graph $G$ (always finite
with unoriented edges and no loops or multiple edges) is a matrix 
obtained by adding an arbitrary diagonal matrix to the adjacency matrix of $G$.

Our first result counts invertible 
Schr\"odinger operators over finite fields for trees (graphs without closed
non-trivial paths):

\begin{thm}\label{mainthma} The number of invertible Schr\"odinger operators
of a finite tree $T$ with $n$ vertices over the finite field 
$\mathbb F_q$ is given by
\begin{eqnarray}\label{formulaSTcharpoly}
&\left(-\sqrt{-q}\right)^n\chi_T\left(\sqrt{-q}+1/\sqrt{-q}\right)&
\end{eqnarray}
where $\chi_T=\det(x\hbox{Id}_n-A)\in\mathbb Z[x]$ 
is the characteristic polynomial of the adjacency matrix $A$ of $T$.
\end{thm}

Theorem \ref{mainthma} is wrong for arbitrary graphs: It fails to
yield integral evaluations at prime-powers for non-bipartite graphs.
It is also wrong for bipartite graphs: Formula (\ref{formulaSTcharpoly})
amounts to $q^4-2q^2+1$
for the the $4$-cycle $C_4$
which has $(q-1)(q^3-q-1)=q^4-q^3-q^2+1$ invertible Schr\"odinger 
operators over $\mathbb F_q$ by Theorem \ref{thmcycles}.

The key-ingredient for proving Theorem \ref{mainthma} is the notion 
of local invariants, a general framework for 
computing invariants of finite (plane) trees.

Our next result enumerates invertible 
Schr\"odinger operators for the $n$-cycle $C_n$ defined as the unique connected graph consisting of $n\geq 3$ vertices of degree $2$:

\begin{thm}\label{thmcycles}
The number $S_{C_n}$ 
of invertible Schr\"odinger operators for the $n$-cycle $C_n$ 
over $\mathbb F_q$ is given by
\begin{eqnarray*}
S_{C_{2n+1}}&=&q^{2n+1}-\frac{1-q^{2n+2}}{1-q^2},\\
S_{C_{4n}}&=&q^{4n}-q^{2n}+\frac{(1-q^{2n})(1-q^{2n+1})}{1-q^2},\\
S_{C_{4n+2}}&=&q^{4n+2}+\frac{(1-q^{2n+1})(1-q^{2n+2})}{1-q^2}.\\
\end{eqnarray*}
if $q$ is odd and by 
\begin{eqnarray*}
S_{C_{2n+1}}&=&q^{2n+1}-\frac{1-q^{2n+2}}{1-q^2},\\
S_{C_{2n}}&=&q^{2n}-q^{n}+\frac{(1-q^{n})(1-q^{n+1})}{1-q^2},\\
\end{eqnarray*}
if $q$ is even.
\end{thm}

Observe that $S_{C_n}$ is  polynomial in $q$, except if $n\equiv 2\pmod 4$
where it is given by two polynomials, depending on the parity of $q$.

The proof of Theorem \ref{thmcycles} is essentially an identity in 
$\mathbb Z[\mathrm{SL}_2(\mathbb F_q)]$, see Theorem \ref{thmgroup}. 
This identity is of independent interest: It 
yields for example Bose-Mesner algebras and 
a good generator of random elements in 
$\mathrm{SL}_2(\mathbb F_q)$. 

Invertible Schr\"odinger operators for the complete graph $K_n$ 
on $n$ vertices
are invertible matrices of size $n\times n$ with arbitrary 
diagonal coefficients and with all off-diagonal coefficients equal to $1$.
The following result gives their number over finite fields:

\begin{thm}\label{thmcpltegraph}
The number of invertible Schr\"odinger operators over $\mathbb F_q$
associated to 
the complete graph on $n$ vertices is given by
$$\frac{(q-1)^{n+1}+(-1)^n}{q}+n(q-1)^{n-1}.$$
\end{thm}

The content of the paper is organized as follows:

Section \ref{sectlocalinvariants}
introduces and gives examples of local invariants, 
the main tool for proving 
Theorem \ref{mainthma}, established in 
Section \ref{sectmainexple}. 

Section \ref{sectaddprop} describes a few additional properties of the 
polynomial $S_T(q)$ enumerating invertible Schr\"odinger operators 
over $\mathbb F_q$ of a finite tree $T$.

Section \ref{secttreecolvert} extends local invariants to trees 
having coloured vertices.

Section \ref{sectionedgesubdivision} studies the behaviour of the polynomial
$S_T(q)$ (defined by Theorem  \ref{mainthma}) under edge-subdivisions.

Section \ref{sectjacobi} refines $S_T$ in order to count invertible 
Schr\"odinger operators (of trees) over finite fields according 
to multiplicities of values of the Jacobi symbol on the diagonal.

In Section \ref{sectgroup}
we give formulae for the coefficients of
\begin{eqnarray}
\left(
\sum_{\mu\in M}\sum_{x\in\mathbb F_q}\left[\left(\begin{array}
{cc}x&\mu\\
-1/\mu&0\end{array}\right)\right]\right)^n\in\mathbb Z[\mathrm{SL}_2(\mathbb F_q)]
\end{eqnarray}
where $M$ is a subgroup of the multiplicative group $\mathbb F_q^*$ of units
in $\mathbb F_q$. Theorem \ref{thmcycles} is a rather straightforward 
consequence of these formulae, as shown in Section \ref{sectproofthmcycles}.

Section \ref{sectproofcompletegraphs} gives an easy proof of 
Theorem \ref{thmcpltegraph}.

A short last Section \ref{sectcomplements} contains a few final remarks.

\section{Local invariants}\label{sectlocalinvariants}

\subsection{Local construction of trees}\label{sectlocalconstrtrees}

We denote by $\mathcal T$ the set of all finite trees and by 
$\mathcal R$ the set of all finite rooted trees. 
%Unless stated otherwise, (rooted) trees will be finite in the sequel. 

Every element of $\mathcal R$ can be constructed (generally not uniquely)
in a finite number of steps involving the following operations:

\begin{description} 
\item[($V$)] Creating a trivial rooted tree consisting of a unique root-vertex.
\item[($E$)] Extending a rooted tree by gluing one end of
an additional edge to the root-vertex and by moving the root vertex to 
the other end of the newly attached edge.
\item[($M$)] Merging two rooted trees by gluing their root-vertices into the 
root-vertex of the resulting tree.
\end{description}

The operation $E$ increases the number of edges and vertices by $1$. 
The operation $M$, applied to two rooted trees having respectively 
$a$ and $b$ vertices
produces a rooted tree with $a+b-1$ vertices and $a+b-2$ edges.

$V$ is constant (independent of any arguments), $E$ operates on elements of $\mathcal R$.
The map $M$ defines a commutative and associative product 
which turns $\mathcal R$ into a commutative monoid with identity $V$
representing the trivial rooted tree reduced to the root vertex.
The monoid $(\mathcal R,M)$ is $\mathbb N$-graded:
The degree of a rooted tree $R$ is the number of non-root vertices in $R$.
The sum over all possible contractions of an edge starting at the root vertex
defines a derivation of degree $-1$ on the graded
monoid-ring $\mathbb K[\mathcal R]$ (over a commutative ring or
field $\mathbb K$).
The map $E$ can thus be thought of as an \lq\lq integral operator\rq\rq
on $\mathcal R$. Algebraically, $\mathbb K[\mathcal R]$ is the free 
commutative algebra with generators $\{E(R)\}_{R\in \mathcal R}$.
Its Hilbert series $\sum_{n=0}^\infty \alpha_nx^n$ encoding the dimension
$\alpha_n$ of homogeneous polynomials of degree $n$ in $\mathbb K[\mathcal R]$
satisfies the identity
$$\sum_{n=0}^\infty \alpha_nx^n=\prod_{n=0}^\infty \left(\frac{1}{1-x^{n+1}}
\right)^{\alpha_n},$$
appearing already in \cite{Ca},
and starts as $1,1,2,4,9,20,48,115,286$, see sequence A81 of \cite{OEIS}.

Finally, the ``forget'' operator, 
\begin{description} 
\item[($F$)] Forgetting the root-structure by turning the root of a rooted tree into an ordinary vertex,
\end{description}
induces a surjection from $\mathcal R$ onto $\mathcal T$.
We have the identity
\begin{eqnarray}\label{Fidentity}
F(M(A,E(B)))=F(M(E(A),B))
\end{eqnarray}
(for all $A,B$ in $\mathcal R$) 
mirroring the fact that an ordinary tree with $n$ vertices 
can be rooted at $n$ different vertices. Identity (\ref{Fidentity})
amounts to requiring $(A,B)\longmapsto F(M(A,E(B)))$ to be symmetric
in its arguments $A$ and $B$.
\begin{figure}[h]
\center{\includegraphics[width=.6\linewidth]{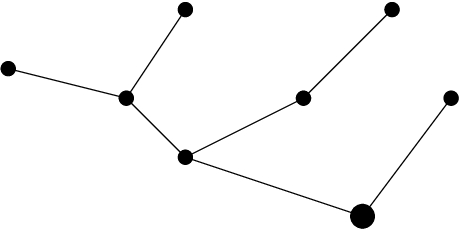}}
\caption{An example of a tree, rooted at the largest dot.}\label{figure1}
\end{figure}

Figure \ref{figure1} 
shows the rooted tree encoded (for example) by
$$M\lbrace E(M\lbrace E(M\lbrace E(V),E(V)\rbrace),
E(E(V))\rbrace ),E(V)\rbrace$$
with curly brackets enclosing arguments of $M$.

\subsection{Digression: plane trees}

A \emph{plane tree} is a tree embedded in the oriented plane, up to orientation-preserving homeomorphisms. Plane trees are abstract trees together with cyclic orders
on sets of edges sharing a common vertex.

A \emph{rooted plane tree} is a plane tree having a root 
together with a refinement into a linear order of the cyclic order 
on root-edges containing the root-vertex.

Rooted plane trees can be constructed using the operators $V,E,M$ already 
considered, except that the associative product $M$ is no longer commutative. 
The operator $F$ turning the root-vertex
into an ordinary vertex
satisfies (\ref{Fidentity}) and the \lq\lq trace-identity\rq\rq
\begin{eqnarray}\label{traceid}
F(M(A,B))=F(M(B,A)).
\end{eqnarray}

The set $\Pi$ of all finite rooted plane trees is a non-commutative
monoid. It is again graded (by the number of non-root vertices) 
and $\mathbb K[\Pi]$ is a non-commutative differential algebra.
Algebraically, $\mathbb K[\Pi]$ is the free 
non-commutative algebra with generators $\{E(R)\}_{R\in \Pi}$. Its Hilbert
series is the algebraic function $C=\sum_{n=0}^\infty c_nx^n=\frac{1}{1-xC}$
whose coefficients define the famous sequence $1,1,2,5,14,\dots,c_n={2n\choose n}\frac{1}{n+1}$ of Catalan numbers, see sequence A108 of \cite{OEIS}.

\subsection{Local invariants of (rooted) trees}

A \emph{local invariant of rooted trees} with values in a 
commutative monoid $\mathbb E$
is a map $i:\mathcal R\longrightarrow \mathbb E$ which can be computed by replacing the construction operators $V,E,M$ by
$v,e,m$ where $v=1$ is the multiplicative identity of $\mathbb E$,
where $e:\mathbb E\longrightarrow \mathbb E$ is an arbitrary map and where 
$m:\mathbb E\times \mathbb E\longrightarrow  \mathbb E$ is the 
product of the monoid $\mathbb E$.

A \emph{local invariant of trees} is a map $f\circ i$ from
$\mathcal T$ to a set of values $\mathbb F$
where $i$ is a local invariant of rooted trees given 
by maps $v,e,m:\mathbb E^*\longrightarrow \mathbb E$ 
as above (with $\mathbb E^*$ denoting respectively
$\emptyset,\mathbb E$ and $\mathbb E^2$)
and where $f:\mathbb E\longrightarrow \mathbb F$ satisfies the
identity
\begin{eqnarray}\label{fidentity}
f(m(A,e(B)))=f(m(e(A),B))
\end{eqnarray}
corresponding to (\ref{Fidentity})
for all $A,B$ in $\mathbb E$.

A trivial example with $\mathbb E=\mathbb F=\mathbb N$ is given by 
$v=0$, $e(x)=x+1$, $m(x,y)=x+y$ and $f(x)=x$.
It counts the number of edges (given by $n-1$ for a tree with $n$ vertices)
of a tree. Replacing $f$ with $f_1(x)=x+1$ we count 
vertices instead of edges.

\begin{rem} The definition of local invariants for (rooted) trees is 
tautological: Every map $A:\mathcal T\longrightarrow \mathbb F$
is a local invariant on the set $\mathcal T$ of all finite trees
by taking $\mathbb E=\mathcal R$ and 
$v=V,e=E,m=M,f=A\circ F$. We are of course interested in local invariants
where the maps $v,e,m$ and $f$ are simple, e.g. given by algebraic
operations over some commutative monoid $\mathbb E$ with a rich algebraic 
structure. 

The terminology \lq\lq local\rq\rq\ alludes to the fact that 
local invariants can be computed algorithmically using \lq\lq local\rq\rq\
operations which modify only neighbourhoods of root-vertices.
\end{rem}

\subsection{Examples of local invariants}\label{subsectclassexples}

\subsubsection{Enumeration of (maximal) 
independent sets} A subset $\mathcal I$ 
of vertices in a graph $G$ is \emph{independent} if two distinct elements
of $\mathcal I$ are never adjacent.
The polynomial $\sum_j\alpha_jx^j$ encoding the number $\alpha_j$
of independent sets with $j$ vertices in a finite tree
can be computed as a local invariant using
\begin{eqnarray*}
v&=&(1,x),\\
e(a,b)&=&(a+b,xa),\\
m((a,b),(\alpha,\beta))&=&(a\alpha,\frac{1}{x}b\beta),\\
f(a,b)&=&a+b.
\end{eqnarray*}
We leave the easy details to the reader. (Hint: The first
coefficient $a$ of $(a,b)$ counts independent sets without the root 
of rooted trees, the second coefficient $b$ counts independent sets 
containing the root-vertex.)

For the tree underlying Figure \ref{figure1} we get
$$1+8x+21x^2+22x^3+8x^4+x^5\ .$$

An independent set $\mathcal I$ of a graph $G$ is \emph{maximal}
if every vertex of $G$ is at most at distance $1$ to $\mathcal I$
(i.e. a vertex $v$ is either in $\mathcal I$ or adjacent to an element
of $\mathcal I$). 

The polynomial $\sum_j\beta_jx^j$ encoding the number $\beta_j$
of maximal independent sets with $j$ vertices in a finite tree
can be computed as a local invariant using
\begin{eqnarray*}
v&=&(1,0,x),\\
e(a,b,c)&=&(b,c,x(a+b)),\\
m((a,b,c),(\alpha,\beta,\gamma))&=&(a\alpha,a\beta+b\alpha+b\beta,
\frac{1}{x}c\gamma),\\
f(a,b,c)&=&b+c.
\end{eqnarray*}
($a$ encodes non-maximal independent sets $\mathcal I$ not containing the 
root $r$ of a rooted tree $R$ such that $\mathcal I\cup\{r\}$ is
maximal independent in $R$ and $\mathcal I$ is maximal independent in the forest 
$R\setminus \{r\}$, the coefficient $b$ encodes maximal independent sets $\mathcal I$ of $R$
such that $r\not\in \mathcal I$ and $c$ encodes maximal independent
sets of $R$ containing the root vertex $r$).

For the tree underlying Figure \ref{figure1} we get
$$4x^3+3x^4+x^5\ .$$

\subsubsection{Enumeration of (maximal) matchings} \label{sssectmatchings}
A \emph{matching} of a graph
is a set of disjoint edges. The polynomial $\sum_j \alpha_j x^j$
with $\alpha_j$ counting
matchings involving $j$ edges can be computed as the local invariant
\begin{eqnarray*}
v&=&(1,0),\\
e(a,b)&=&(a+b,xa),\\
m((a,b),(\alpha,\beta))&=&(a\alpha,a\beta+b\alpha),\\
f(a,b)&=&a+b.
\end{eqnarray*}
We leave the easy details to the reader. (Hint: The first
coefficient $a$ of $(a,b)$ counts matchings of a rooted tree 
not involving the root, the second coefficient $b$ counts matchings
involving the root.)

For the tree underlying Figure \ref{figure1} we get
$$1+7x+13x^2+7x^3\ .$$

A matching of a graph $G$ is \emph{maximal} if it intersects every edge of $G$.
The polynomial $\sum_j \beta_j x^j$
with $\beta_j$ counting maximal
matchings involving $j$ edges can be computed as the local invariant
\begin{eqnarray*}
v&=&(0,1,0),\\
e(a,b,c)&=&(b,c,x(a+b)),\\
m((a,b,c),(\alpha,\beta,\gamma))&=&(a\alpha+a\beta+b\alpha,b\beta,(a+b)\gamma+
c(\alpha+\beta)),\\
f(a,b,c)&=&b+c
\end{eqnarray*}
($a$ counts not maximal matchings of a rooted tree $R$ inducing maximal
matchings on the forest $R\setminus \{r\}$ obtained by removing the root 
$r$ from $R$, the coefficient $b$ counts maximal matchings of $R$ not
involving the root and $c$ counts maximal matchings of $R$ involving the root $r$).

For the tree underlying Figure \ref{figure1} we get
$$7x^3\ .$$

\subsubsection{The characteristic polynomial of the adjacency matrix}
\label{sssectcharpoly}

We write $(a,b)\in\mathbb Z[x]$ if the characteristic polynomial of a rooted
tree is given by $ax-b$ with $ax$ corresponding to the contribution of 
the diagonal entry
associated to the root. Elementary matrix-transformations
show that the characteristic polynomial
$\det(x\hbox{ Id}-A)$ of the adjacency matrix $A$ 
is a local invariant defined by 
\begin{eqnarray*}
v&=&(1,0),\\
e(a,b)&=&(xa-b,a),\\
m((a,b),(\alpha,\beta))&=&(a\alpha,a\beta+b\alpha),\\
f(a,b)&=&xa-b.
\end{eqnarray*}
For the tree underlying Figure \ref{figure1} we get
$$x^8-7x^6+13x^4-7x^2\ .$$

Similarly,
\begin{eqnarray*}
v&=&(1,0,0),\\
e(a,b,c)&=&(-(x+c+1)a+b,-a,1),\\
m((a,b,c),(\alpha,\beta,\gamma))&=&(a\alpha,a\beta+b\alpha,c+\gamma),\\
f(a,b,c)&=&-(x+c)a+b
\end{eqnarray*}
computes
the characteristic polynomial of the combinatorial Laplacian (given
by $D_{\mathrm{deg}}-A$ where $A$ is the adjacency matrix and $D_{\mathrm{deg}}$
is the diagonal matrix defined by vertex-degrees) of 
a tree.

For the tree underlying Figure \ref{figure1} we get
$$x^8+14x^7+76x^6+204x^5+286x^4+204x^3+67x^2+8x\ .$$

\subsection{An example with values in $\mathbb N[z_0,z_1,\dots]$}\label{subsectimprtntexple}

Given a formal power series $Z=\sum_{i=0}^\infty z_it^i\in \mathbb A[[t]]$
and a polynomial $B=\sum_{i=0}^N b_it^i\in \mathbb A[t]$
with coefficients in a commutative ring $\mathbb A$, we denote
by 
$$\langle Z,B\rangle=\sum_{i=0}^N z_ib_i$$
the scalar product of the \lq\lq coefficient-vectors\rq\rq
$(z_0,z_1,\dots)$ and $(b_0,\dots,b_N,0,\dots)$.

\begin{prop}\label{propgenlclinv} We denote by $Z=\sum_{i=0}^\infty z_it^i$ a formal power series 
in $t$ with coefficients $z_i$.
The formulae
\begin{eqnarray*}
v&=&1,\\
e(A)&=&\langle Z,A\rangle+tA,\\
m(A,B)&=&AB,\\
f(A)&=&\langle Z,A\rangle
\end{eqnarray*}
define a local invariant of trees with values in 
$\mathbb N[z_0,z_1,z_2,\dots]$, respectively of rooted trees 
with values in $\mathbb N[z_0,z_1,z_2,\dots][t]$.
\end{prop}

For the tree underlying Figure \ref{figure1} we get
\begin{eqnarray*}
&&z_0^8+7z_0^6z_1+8z_0^5z_2+z_0^4(13z_1^2+10z_3)+z_0^3(18z_1z_2+11z_4)\\
&&+z_0^2(7z_1^3+12z_1z_3+3z_2^2+8z_5)+z_0(9z_1^2z_2+8z_1z_4+2z_2z_3+4z_6)\\
&&+z_1^2z_3+2z_1z_2^2+2z_1z_5+z_2z_4+z_7\ .
\end{eqnarray*}
The specialization $z_2=z_3=\dots=0$ of Proposition \ref{propgenlclinv} is
particularly interesting in the sense that it gives an invariant in
$\mathbb N[z_0,z_1]([t])$  of (rooted)
trees which behaves naturally with respect to the differential structure of 
$\mathbb N[\mathcal R]$. Given a rooted tree $R$ 
with invariant $c_R(0)+c_R(1)t+\dots$, the constant coefficient $c_R(0)$
corresponds to $R$ and the linear coefficient $c_R(1)$ corresponds 
to the derivative (as in Section \ref{sectlocalconstrtrees}) of $R$.
The result (with $Z=z_0+z_1t$) for the tree underlying Figure \ref{figure1} 
is
$$z_0^2(z_0^2+z_1)(z_0^4+6z_0^2z_1+7z_1^2)\ .$$

 \noindent{\bf Proof of Proposition \ref{propgenlclinv}} The set $\mathbb N[z_0,z_1,z_2,\dots][t]$
of all polynomials is a multiplicative monoid with product 
$m(A,B)=AB$ and identity $1$. The operator $e$ defines a map from 
$\mathbb N[z_0,z_1,z_2,\dots][t]$ into itself. This shows that
the formulae of Proposition \ref{propgenlclinv} define a 
local invariant of rooted trees.

Symmetry in $A,B$ of
\begin{eqnarray*}
f(m(A,e(B)))&=&f(A(\langle Z,B\rangle+tB))\\
&=&\langle Z,A\rangle\langle Z,B\rangle+\langle Z,tAB\rangle
\end{eqnarray*}
implies that identity (\ref{fidentity}) holds.\hfill$\Box$

\subsubsection{Examples} Up to a power of $x$, the example 
of Section \ref{sssectmatchings} enumerating matchings corresponds to 
$Z=\frac{1+t}{x}$ with coefficients
$z_0=z_1=1/x$ and $z_2=z_3=\dots=0$.

Our next example, stated as a theorem, will be the crucial ingredient 
for proving Theorem \ref{mainthma}:

\begin{thm}\label{thmcharpoly}
The characteristic polynomial of Section \ref{sssectcharpoly}
corresponds to $Z=x-t$ with coefficients
$z_0=x,z_1=-1$ and $z_2=z_3=\dots=0$.
\end{thm}

\noindent{\bf Proof} Follows easily from the formulae given in 
Section \ref{sssectcharpoly}.\hfill$\Box$

\section{A local invariant enumerating 
Schr\"odinger operators}\label{sectmainexple}

We consider the local invariant 
$S:\mathcal T\longrightarrow \mathbb Z[q]$ of
trees given by the specialization
$Z=q-1+qt$ with coefficients
$z_0=q-1,z_1=q$ and $z_2=z_3=\dots=0$ of the local invariant described by
Proposition \ref{propgenlclinv}.

Since $S$ depends only on the constant and on the linear coefficient of 
the corresponding local invariant $a+bt+\dots$ of rooted trees, 
we can also define $S$ by the formulae
\begin{eqnarray}\label{formulaeS}
\begin{array}{rcl}
v&=&(1,0),\\
e(a,b)&=&((q-1)a+qb,a),\\
m((a,b),(\alpha,\beta))&=&(a\alpha,a\beta+b\alpha),\\
f(a,b)&=&(q-1)a+qb
\end{array}
\end{eqnarray}
with $(a,b)$ representing the series expansion $a+bt+O(t^2)$.

For the (unrooted) tree underlying Figure \ref{figure1} we get
the polynomial
$$(q-1)^2(q^2-q+1)(q^4+2q^3+q^2+2q+1) \ .$$

\begin{rem} As a mnemotechnical device, the formula for $m$ corresponds also to 
the addition $\frac{b}{a}+\frac{\beta}{\alpha}=\frac{a\beta+b\alpha}{a\alpha}$
with forbidden simplification and the formula for $e$ is, up to simplification 
by $a$, given by the homography $\frac{b}{a}\longmapsto
\left(\begin{array}{cc}0&1\\q&q-1\end{array}\right)\frac{b}{a}=\frac{1}
{q\frac{b}{a}+q-1}$. The formula for $f$ can be recovered from $e$ using the 
identity $e(a,b)=(f(a,b),a)$.
\end{rem}

We write $S_T$ for the local invariant in $\mathbb Z[q]$ associated to a tree
$T$. Similarly, given a rooted tree $R$, 
we denote by $S_R=(a,b)\in\left(\mathbb Z[q]\right)^2$
the corresponding pair of polynomials defined by formulae (\ref{formulaeS}).

\begin{prop}\label{propSlocal} $S_T$ counts the number of invertible Schr\"odinger operators for a finite tree $T$
over the finite field $\mathbb F_q$.
\end{prop}

A matrix $M$ with rows and columns indexed by vertices of a
graph $G$ is a \emph{$G$-matrix}
if non-zero off-diagonal coefficients $m_{s,t}$ of $M$ correspond to edges
$\{s,t\}$ of $G$.
Diagonal entries of $G$-matrices are arbitrary.
The off-diagonal support (set of non-zero coefficients) 
of a $G$-matrix encodes thus the edge-set of $G$.
A $G$-matrix of an unoriented graph $G$ 
has always a symmetric support but is not necessarily symmetric.
We have:

\begin{prop}\label{propinvTmatrix} The number of invertible $T$-matrices over 
$\mathbb F_q$ of a finite tree $T$ having
$n$ vertices is given by $(q-1)^{2n-2}S_T$.
\end{prop}

\noindent{\bf Proof of Proposition \ref{propSlocal}}
We consider the obvious action
on $T$-matrices of the abelian group $\left(\mathbb F_q^*\right)^n\times \left(\mathbb F_q^*\right)^n$ of pairs of invertible diagonal matrices (with coefficients in 
$\mathbb F_q$) by left and right multiplication. An orbit of a $T$-matrix 
without non-zero diagonal coefficients contains exactly
one Schr\"odinger operator stabilized by a subgroup of order $(q-1)^2$.
All other orbits contain exactly $q-1$ different Schr\"odinger operators, 
each stabilized by a subgroup of order $q-1$. The number of (invertible) 
$T$-matrices is thus exactly $(q-1)^{2n-2}$ times larger than 
the number of (invertible) Schr\"odinger operators.
The result follows now from Proposition \ref{propinvTmatrix}.
\hfill$\Box$

\noindent{\bf Proof of Proposition \ref{propinvTmatrix}} 
An \emph{$R$-matrix} for a rooted tree $R\in\mathcal R$ 
is a $T$-matrix for the underlying
unrooted tree $T$ with an unknown $x$ on the diagonal 
corresponding to the root of $R$.
The determinant of an $R$-matrix over a finite field $\mathbb F_q$ 
is an affine function of the form 
$ax+b\in\mathbb F_q[x]$. 

We consider now a fixed rooted tree $R$.
Given two subsets $\mathcal A,\mathcal B$ of $\mathbb F_q$, we denote by 
$\nu(\mathcal A,\mathcal B)$ the number
of $R$-matrices of determinant $ax+b$ with 
$(a,b)\in \mathcal A\times \mathcal B$.
We encode the natural integers $\nu(0,0),\nu(0,\mathbb F_q^*),\nu(\mathbb F_q^*,0),
\nu(\mathbb F_q^*,\mathbb F_q^*)$ (with $0$ denoting the singleton 
subset $\{0\}$ consisting of the zero-element in $\mathbb F_q$) using the square matrix
$$\left(\begin{array}{cc}
 \nu(0,0)&\nu(0,\mathbb F_q^*)\\
\nu(\mathbb F_q^*,0)&
\nu(\mathbb F_q^*,\mathbb F_q^*)
\end{array}\right).$$
Right and left multiplications by invertible diagonal matrices 
essentially preserve 
the set of $R$-matrices. More precisely, this holds up to replacement
of the unknown $x$ by a non-zero multiple $\lambda x$ (with $\lambda\in
\mathbb F_q^*$) of it.
Since $x$ can be thought of as a simple place-holder for an
arbitrary element of $\mathbb F_q$, such a scalar $\lambda$ can be dismissed.
It follows that we have $\nu(\lambda a,\mu b)=\nu(a,b)$
(using a slight abuse of notation)
if $\lambda$ and $\mu$  belong both to the set $\mathbb F_q^*$
of invertible elements in $\mathbb F_q$.
Elementary linear algebra shows now that the operators $V,E,M,F$ correspond to
the operators
\begin{eqnarray*}
v&=&\left(\begin{array}{cc}0&0\\1&0\end{array}\right),\\
e&=&\left(\begin{array}{cc}a&b\\c&d\end{array}\right)\longmapsto
(q-1)^2\left(\begin{array}{cc}qa&c+d\\qb&(q-1)(c+d)\end{array}\right),\\
m&=&(\left(\begin{array}{cc}a&b\\c&d\end{array}\right),
\left(\begin{array}{cc}\alpha&\beta\\\gamma&\delta\end{array}\right))\longmapsto
\left(\begin{array}{cc}A&B\\C&D\end{array}\right),\\
f&=&\left(\begin{array}{cc}a&b\\c&d\end{array}\right)\longmapsto
qb+(q-1)(c+d)
\end{eqnarray*}
where 
\begin{eqnarray*}
A&=&a\alpha+a\beta+b\alpha+a\gamma+c\alpha+\\
&&+a\delta+d\alpha+b\beta,\\
B&=&b\gamma+c\beta+b\delta+d\beta,\\
C&=&c\gamma+\frac{1}{q-1}d\delta,\\
D&=&c\delta+d\gamma+\frac{q-2}{q-1}d\delta.
\end{eqnarray*}
The coefficient $\nu(0,0)$ contributes never to the number 
of invertible Schr\"odinger operators and can be discarded.
Inspection of the above formulae shows that we can lump together 
$\nu(\mathbb F_q^*,0)$ and $\nu(\mathbb F_q^*,\mathbb F_q^*)$
into a first coordinate with the second coordinate given by 
$\nu(0,\mathbb F_q^*)$.
This leads to the formulae (\ref{formulaeS}) 
for $S_T$, except for an extra factor
of $(q-1)^2$ for every edge of $T$.\hfill$\Box$

\subsection{Proof of Theorem \ref{mainthma}}

Substituting $x$ with $\sqrt{-q}+1/\sqrt{-q}$ and multiplying by 
the correct sign and power of $\sqrt{-q}$, the
expression $(-\sqrt{-q})^n\chi_T\left(\sqrt{-q}+1/\sqrt{-q}\right)$
(with $n$ denoting the number of vertices of $T$)
can be computed using Theorem \ref{thmcharpoly} as the local invariant
given by 
\begin{eqnarray*}
v&=&1,\\
e(a+tb+O(t^2))&=&-\sqrt{-q}\left(\sqrt{-q}+1/\sqrt{-q}\right)a-\sqrt{-q}^2b
+at+O(t^2)\\
&=&(q-1)a+qb+at+O(t^2),\\
m(A,B)&=&AB,\\
f(a+tb+O(t^2))&=&-\sqrt{-q}\left(\sqrt{-q}+1/\sqrt{-q}\right)a-\sqrt{-q}^2b
\\
&=&(q-1)a+qb\\
\end{eqnarray*}
and coincides thus with the local invariant given by formulae
(\ref{formulaeS}) which define the counting function $S_T$ for Schr\"odinger
operators by Proposition \ref{propSlocal}.\hfill$\Box$

\section{A few additional properties of $S_T$}\label{sectaddprop}

Since Theorem \ref{mainthma} links $S_T$ closely to the characteristic polynomial (of an 
adjacency matrix), many properties of $S_T$ mirror
properties of characteristic polynomials for trees.
For example, since characteristic polynomials (of integral 
matrices) are monic and integral, the polynomial $S_T$ is monic and
integral.

A property not linked to the characteristic polynomial but
to the formula (\ref{formulaSTcharpoly}) defining $S_T$
is the fact that the polynomial $S_T$ associated to a tree
$T$ with $n$ vertices satisfies the equation
\begin{eqnarray}\label{signdegrepal}
S_T(q)&=&(-q)^nS_T(1/q).
\end{eqnarray}
We call this property \emph{sign-degree-palindromicity}.
It implies that a complex number $\rho$ is a root of $S_T$
if and only if $1/\rho$ is a root.

Another easy fact (left to the reader), is the observation that $S_T$ is
always of the form $q^n-q^{n-1}+\dots$. This is (up to $O(q^{n-2})$)
the expected number of non-zero elements among $q^n$ (uniformly distributed)
random elements of $\mathbb F_q$.

\subsection{Root-locus of $S_T$}\label{subsectrootlocus}

All roots of the characteristic polynomial of a graph are real.
Since a tree $T$ is bipartite, a real number $\rho$ is a root of
$\chi$ (the characteristic polynomial of an adjacency matrix of $T$)
if and only if $-\rho$ is a root. 
A pair of non-zero roots $\pm \rho$ of $\chi$ gives thus rise
to the pair $\sigma,1/\sigma$ of roots of $S_T$ satisfying
the equation
\begin{eqnarray*}
0&=&(\sqrt{-\sigma}+1/\sqrt{-\sigma}-\rho)(\sqrt{-\sigma}
+1/\sqrt{-\sigma}+\rho)\\
&=&-\sigma+2-1/\sigma-\rho^2.
\end{eqnarray*}
We have thus
$$\sigma^{\pm 1}=\frac{2-\rho^2\pm \sqrt{(\rho^2-2)^2-4}}{2}.$$
For $\rho\in [-2,2]$ we get two conjugate roots $\sigma,\overline{\sigma}=1/\sigma$ on the complex unit circle (except perhaps for $\rho=0$ giving sometimes
rise to a unique root $1$ of $S_T$), for $\rho$ of absolute value larger than 
$2$ we get two negative real roots $\sigma,1/\sigma$.

All roots of $S_T$ are thus on the union of the complex unit circle
with the real negative half-line.

Trees with all roots of $S_T$ on the unit circle are subtrees of affine 
Dynkin diagrams of type $D$ or $E$, see for example \cite{MS05}
which describes also all trees such that $S_T$ has exactly one real root
$<-1$. (More precisely, \cite{MS05} deals with
polynomials of the form (\ref{formulaSTcharpoly}) (up to trivial signs)
which give rise to Salem numbers.)

Our next result shows that ``simple''\ trees give rise to polynomials
$S_T$ with few real negative roots:

\begin{prop}\label{cornumbnegrootsab}
The number of real negative zeroes of $S_T\in\mathbb Z[q]$
associated to a tree $T$ 
is at most equal to twice the number of non-root vertices of degree at least 
$3$ in $T$.
\end{prop}

The proof uses the following auxiliary result:

\begin{lem}\label{propnumbnegrootsab} The numbers of real negative zeroes of the polynomials
$a,b$  associated by $S_R=(a,b)$ to a rooted tree $R$ 
are at most equal to twice the number of non-root vertices of degree at least 
$3$ in $R$.
\end{lem}

\noindent{\bf Proof of Proposition \ref{cornumbnegrootsab}}
We turn $T$ into a rooted tree $R$ by choosing a root vertex 
at a leaf of $T$. The result follows now by applying 
Lemma \ref{propnumbnegrootsab} to the first polynomial
of $S_{E(R)}=(S_T,*)$.\hfill$\Box$

\noindent{\bf Proof of Lemma \ref{propnumbnegrootsab}}
Let $R$ be a rooted tree. If the root vertex $v_*$ of $R$ is not
a leaf, then $S(R)=(A,B)=(a\alpha,a\beta+b\alpha)$ where $(a,b)$ and 
$(\alpha,\beta)$ are associated to 
smaller non-trivial rooted trees $R_1,R_2$ such that 
$R=M(R_1,R_2)$. The result holds thus for $A=a\alpha$ 
by induction on the number of vertices
and it holds for $B$ since roots of $A$ and $B$ interlace in an obvious sense
on $\mathbb S^1\cup \mathbb R_{<0}$.
If the root vertex $v_*$ is a leaf, the result holds by a 
straightforward computation if $R$ 
is a leaf-rooted path (Dynkin diagram of type $A$).
Otherwise, the tree $R$ contains a vertex $w$ of degree at least $3$.
Working with the rooted tree $R_w$ corresponding to $T$ rooted at $w$,
we see that $a_w,b_w$ with $S(R_w)=(a_w,b_w)$ have at most $2(k-1)$
real negative zeroes where $k$ is the number of vertices of degree
at least $3$ in $T$. This implies that $S_T$ has at most $2+2(k-1)=2k$
real negative zeroes. Since $S_T$ has at least as many real 
negative zeroes as $b$ involved in $S(R)=(a,b)$ and since roots of $a$
and $b$ interlace, the polynomials $a,b$ have both 
at most $2k$ real zeroes. \hfill$\Box$

\subsection{Multiple roots of $S_T$}

A vertex $w$ of degree $k$ of a tree $T$ can be considered as the
result of gluing $k$ maximal leaf-rooted subtrees of $T$
along their root-leaf corresponding to $w$. This 
construction is linked to some multiple roots of $S_T$
as follows: Given a leaf-rooted tree $R$, let $w_1,\dots,w_k$
be the list of all vertices of $T$ involving $R$ (i.e. at least one
of the maximal leaf-rooted subtrees of $T$ with root
$w_i$ is isomorphic to $R$).
Given such a vertex $w_i$, we denote by $r_i+1\geq 1$
the number of occurrences of $R$ at $w_i$.
We have the following result:

\begin{prop}\label{propmult} The polynomial $S_T$ is (at least) divisible 
by $a^r$ where $r=r_1+\cdots+r_k$ and where $a$ is defined by $S_R=(a,b)$.
\end{prop}

\noindent{\bf Proof} Gluing $k_i+1$ copies of $R$  along their
root gives a rooted tree with invariant $(a^{k_i+1},(k_i+1)a^{k_i}b)$.
Linearity of the formulae for $e,f$ and $m$ 
implies now the result.\hfill$\Box$

Proposition \ref{propmult} explains the factor $q^2-q+1$ and one of the 
factors $q-1$ of the polynomial 
$S_T=(q-1)^2(q^2-q+1)(q^4+2q^3+q^2+2q+1)$ with $T$ given by Figure 
\ref{figure1}. Since $S_T$ is palindromic, the factor
 $q-1$ divides $S_T$ with even multiplicity.
All cyclotomic factors in this example have thus easy explanations.

\section{Trees with coloured vertices}\label{secttreecolvert}

A rooted tree with ordinary (non-root) vertices coloured (not necessarily properly,
i.e. adjacent vertices might have identical colours) by a set 
$\mathcal C$ can be constructed using the construction-operators $V$
(creation of a root-vertex), $M$ (merging of two rooted trees 
along their root) and
replacing $E$ by operators $E_c$ (for $c\in\mathcal C$)
depending on the final colour of 
the initial root-vertex. For ordinary trees, one works with operators 
$F_c$ indexed by all possibilities of colouring the root-vertex after turning
it into an ordinary vertex.

Identity (\ref{fidentity}) has to be replaced by
\begin{eqnarray}\label{colouredfidentity}
F_s(M(A,E_t(B))&=&F_t(M(E_s(A),B))
\end{eqnarray}
for all $s,t\in \mathcal C$ and for all $A,B\in \mathcal R$.

Local invariants for coloured (rooted or plane) trees are 
defined in the obvious way.

\subsection{A few examples of coloured invariants}

\subsubsection{Colourings defined by (virtual) rooted trees}\label{sectcolvirttrees}
Every local invariant gives rise to a coloured local invariant by 
chosing colour-constants $u_c\in \mathbb E$ for all colours $c\in \mathcal C$
and by replacing
$e$ with $e_c(A)=e(m(A,u_c))$ and  $f$ with $f_c(A)=f(m(A,u_c))$. 
These invariants
amount to attaching ``virtual trees'' $U_c$ corresponding to $u_c$ and
representing colours to all ordinary vertices. 

A particular case, closely related to $S_T$, will be discussed in Section 
\ref{sectessfinitetrees}.

\subsubsection{A coloured local invariant with values in $\mathbb N[z_0(c),z_1,
z_2,z_3,\dots]$}\label{sectcollocalZ}
The example of Proposition \ref{propgenlclinv}
in Section \ref{subsectimprtntexple} can easily be generalized to 
a coloured local invariant by considering formal power-series
$Z(s)=z_0(s)+\sum_{i=1}^\infty z_it^i$ with constant terms (with respect to 
$t$) 
depending on colours. Identity (\ref{colouredfidentity}), corresponding to
$$\langle Z(s_1),A\rangle\langle Z(s_2),B\rangle+\langle Z(s_1),tAB\rangle
=\langle Z(s_1),A\rangle\langle Z(s_2),B\rangle+\langle Z(s_2),tAB\rangle,$$
holds since $\langle Z(s),tC\rangle$ does not depend on the colour $s$. 

As an example, we can consider the invariant given by $Z(v)=x_v+qt$
generalizing the invariant counting invertible Schr\"odinger 
obtained by the specialization $x_v=q-1$ for all $v$, see 
Section \ref{sectmainexple}.

\subsubsection{Bicoloured characteristic polynomial}

Trees are connected bipartite graphs and have thus a canonical 
proper $2$-colouring or bi-colouring, well-defined up to colour-exchange. 
The corresponding bi-coloured variation of the characteristic polynomial of the adjacency matrix
is given by computing the determinant of the matrix coinciding with the adjacency matrix outside the diagonal and with diagonal coefficients $-x$ or $-y$
according to the bipartite class of the corresponding vertex. The resulting determinant is well-defined in $\mathbb Z[x,y]$ up to exchanging $x$ with $y$ and can be computed 
as a local invariant. This construction works of course also for the combinatorial Laplacian of a tree.

\subsubsection{Coloured Schr\"odinger operators} 

The enumeration of Schr\"odinger operators according to coloured diagonal zeros
leads to a local invariant of coloured trees. It takes
its values in $\mathbb Z[q,\mathcal C]$
with the coefficient (in $\mathbb Z[q]$) of a monomial $\prod_j c_j^{e_j}\in\mathcal C^*$ counting the number of Schr\"odinger operators with $e_j$ zero terms
on diagonal elements associated to vertices of colour $c_j$.

The corresponding operators are given by
\begin{eqnarray*}
&&v=(0,1,0),\\
&&e_s(a,b,c)=(sb+c,(q-1+s)a,(q-1)b+(q-2+s)c),\\
&&m((a,b,c),(\alpha,\beta,\gamma))\\
&&\ =(a\beta+b\alpha+a\gamma+c\alpha,b\beta+\frac{c\gamma}{q-1},
b\gamma+c\beta+\frac{q-2}{q-1}c\gamma)\\
&&f_s(a,b,c)=(q-1+s)a+(q-1)b+(q-2+s)c.
\end{eqnarray*}
(with $(a,b,c)$ standing for $a=N(0,\mathbb F_q^*),b=N(\mathbb F_q^*,0)$ and
$c=N(\mathbb F_q^*,\mathbb F_q^*)$ 
where $N(*,*)$ is as in the proof of Proposition \ref{propinvTmatrix}).

Skeptical readers are invited to check the identity
$$f_s(m((a,b,c),e_t(\alpha,\beta,\gamma)))=
f_t(m(e_s(a,b,c),(\alpha,\beta,\gamma)))$$
corresponding to (\ref{colouredfidentity}).

\subsection{A further generalization: working with coloured monoids}\label{sectcolmonoid}

Rooted vertex-coloured trees with colours at all vertices including the root
do not form a natural monoid. However the subset of rooted trees with 
a given root colour is obviously a monoid (by gluing, as before,
all root-vertices into a root-vertex of the same colour).

A local invariant for such rooted 
coloured trees is given by monoids $\mathcal M_c$ with identities $v_c$ 
and products $m_c$ and by edge-maps 
$e_{c_i,c_f}:\mathcal M_{c_i}\longrightarrow \mathcal M_{c_f}$
(depending on the root-colour $c_i$ of the initial argument-tree 
and on the colour $c_f$ of the final, added root-vertex)
among all coloured monoids.

Adding maps $f_c:\mathcal M_c\longrightarrow \mathbb F$
into some set $\mathbb F$ such that we have for all pairs of colours $c_1,c_2$
the identity 
\begin{eqnarray}\label{formulecolmonoid}
f_{c_1}(m_{c_1}(A_1,e_{c_2,c_1}(A_2)))=f_{c_2}(m_{c_2}(e_{c_1,c_2}(A_1),A_2))
\end{eqnarray}
(with $A_1\in \mathcal M_{c_1},A_2\in \mathcal M_{c_2}$), 
we get an invariant of coloured trees.

Section \ref{sectjacobi} contains an example of this construction.

\section{Edge-subdivisions}\label{sectionedgesubdivision}

\subsection{Properties of $S$}

\begin{prop}\label{proponevertexlimit}
 Let $T_i$ be a sequence of trees obtained by subdividing
all edges around a fixed vertex $w$ of degree $d\geq 3$ of a finite tree 
$T$ into 
a larger and larger number of edges (by insertion of additional vertices 
of degree $2$). Then there exists a sequence of roots $\rho_i$ of $S_{T_i}$
converging to $1-d$.
\end{prop}

The proof follows easily from the discussions in Section 
\ref{sectessfinitetrees}.

Applying Proposition \ref{proponevertexlimit} at all vertices we get:

\begin{cor}\label{corsubdividingallroots} Subdividing all edges
of a fixed finite tree $T$ leads to a sequence of polynomials
with $k$ strictly negative roots converging to $1-d_i$ 
where $d_1,\dots,d_k$ are the degrees of all $k$ vertices of $T$
with degrees $>2$.
\end{cor}

The density of roots on $\mathbb S^1$ under edge-subdivisions can be shown 
to behave as expected:

\begin{prop}
Given a sequence $T_i$ of finite trees obtained
by subdividing some edges of a given fixed tree $T$ into more and 
more sub-edges, the density of roots of 
$S_{T_i}$ on the unit circle $\mathbb S^1$ converges to the Lebesgue measure 
of $\mathbb S^1$. Otherwise 
stated, the proportion of roots of $S_{T_i}$ in a given 
sub-interval $I$ of $\mathbb S^1$ tends to $\frac{1}{2\pi}\mathop{length}(I)$.
\end{prop}

We omit the proof.

\subsection{Essentially finite trees and limits of real negative roots
of $S_T$ under edge-subdivisions}\label{sectessfinitetrees}

A perhaps infinite 
tree is \emph{essentially finite} if it is obtained by subdiving
edges of a finite tree at most countably many times. Essentially finite
trees can be considered as trees with edges weighted by elements of 
$\{1,2,\dots\}\cup\{\infty\}$. Edge-weights encode the
final number of edges after subdivision. Forbidding vertices of degree
$2$ leads to unique representations of this form.

Essentially finite trees with an infinite number of vertices
have no longer an $S$-polynomial.
They define however a finite set $\rho_1,\dots,\rho_k$ 
of $k$ real numbers $<-1$ where $k$ is at most
equal to the number of 
vertices of degree $\geq 3$ in the following way: Approximate such
a tree $T$ by a sequence $T_i$ of finite trees in the obvious way
(by replacing all infinite edge-weights by large finite edge-weights)
and consider the limits (which exist by Theorem \ref{thmlimitroots})
of all real roots $<-1$ of the polynomials
$S_{T_i}$, taking into account multiplicities. The aim of this section
is to compute these numbers and to study a few of their properties.

An essentially finite tree is \emph{arc-connected} if all 
non-leaves are at finite distance. Equivalently, a finite tree
is arc-connected if all its edges with infinite weight contain a 
leaf.
An \emph{arc-connected component} of an essentially finite 
tree is a subtree defined by all vertices at finite distance of a 
non-leaf. An arc-connected essentially finite tree is a finite tree together 
with attachements of finitely many infinite rays (ending at an 
\lq\lq ideal\rq\rq leaf-vertex) at vertices.
An essentially finite arc-connected tree can be encoded by a finite tree
with $\mathbb N$-weighted vertices. Vertex-weights indicate
numbers of attached infinite rays.

Arbitrary essentially finite trees can be decomposed into essentially finite 
arc-connected trees with two arc-connected components intersecting 
at most in a unique \lq\lq ideal midpoint\rq\rq of an infinitely 
subdivised edge. 

The data 
\begin{eqnarray}\label{formulaecolouredS}
\begin{array}{lcr}
v&=&(1,0),\\
e_k(a,b)&=&(q(-ka+b)+(q-1)a,a),\\
m((a,b),(\alpha,\beta))&=&(a\alpha,a\beta+b\alpha),\\
f_k(a,b)&=&q(-ka+b)+(q-1)a
\end{array}
\end{eqnarray}
%\begin{eqnarray}\label{formulaecolouredS}
%\begin{array}{lcr}
%v&=&(0,1),\\
%e_k(a,b)&=&(b,q(a-kb)+(q-1)b),\\
%m((a,b),(\alpha,\beta))&=&(a\beta+b\alpha,b\beta),\\
%f_k(a,b)&=&q(a-kb)+(q-1)b
%\end{array}
%\end{eqnarray}
defines a local invariant for coloured trees with colour-values in
some ring. We leave it to the reader to check that
(\ref{colouredfidentity}) holds. This local 
invariant corresponds to a coloured invariant where 
a vertex $v$ of weight $k$ is decorated by a non-existent ideal 
rooted tree with $S$-polynomial $(1,-k)$, see also 
Section \ref{sectcolvirttrees}.
It can also be considered as the specialization given by $Z(k)=
(q-1-qk)+qt$ of the invariant described in Section \ref{sectcollocalZ}.

We denote by $Q_T\in \mathbb Z[q]$ the polynomial
defined by (\ref{formulaecolouredS}) for 
a finite tree $T$ with $\mathbb Z-$coloured vertices.

For values in $\mathbb N$ encoding numbers of \lq\lq limit rays''
of an arc- connected essentially finite tree, the polynomial $Q_T$
has the inverse limit values $1/\rho_i\in(-1,0)$
(taking multiplicities into account) among its roots.
In particular, the limit-values $\rho_i$ are algebraic
numbers (in fact algebraic integers since $Q_T(0)\in \{\pm 1\}$)
of degrees bounded by the number of vertices in a
$\mathbb N$-vertex-coloured finite tree representation.
Indeed, given a sequence $R_i$ of increasing leaf-rooted paths (Dynkin diagrams of type $A$), the evaluation of  
$S_{R_i}$ at a complex number $\rho$ of norm $<1$
tends to $\pm\frac{1}{1+\rho}(1,-1)$ (with a sign depending on the parity of
the number of vertices).
% by the results of Section \ref{expletypeA}.
Linearity of the maps $e,m$ and $f$ implies now that the limit values 
$\rho_i$ are algebraic numbers.

Given a sequence $T_i$ of finite trees obtained by subdividing increasingly
often a unique edge of an (essentially finite) tree, the argument above 
shows that a real limit-root $\rho$ \lq\lq belongs\rq\rq to one 
of the two arc-connected components of the essentially finite 
limit-tree.

\noindent{\bf Proof of Proposition \ref{proponevertexlimit}} Follows 
from the fact 
that the essentially finite  arc-connected tree $T$ 
represented by an isolated vertex of weight
$d\geq 3$ gives rise to $Q_T=f(1,-d)=(1-d)q-1$.\hfill$\Box$

We have:

\begin{thm}\label{thmlimitroots} All non-limit roots of $Q_T$ have norm $\geq 1$
if $T$ encodes an essentially finite arc-connected 
tree (i.e. if all vertex-colours are
in $\mathbb N$).
\end{thm}

\noindent{\bf Proof} Let $\sigma$ be a root of $Q_T$ in the open complex
unit disc.
The root $\sigma$ can be approximated with arbitrary accuracy 
by a root of $S_{T'}$ where $T'$ is a finite tree approximating
the essentially finite tree $T$. This implies that $\sigma$ has
to be a real negative number, see Section \ref{subsectrootlocus}.\hfill$\Box$

We call the sum of vertex-weights the \emph{weight-degree} of 
a $\mathbb N$-coloured connected essential tree $T$.

Trees of degree $0$ are ordinary trees. Their $S$-polynomials define in some sense \lq\lq generalized Salem 
numbers\rq\rq. Trees of degree $1$ define \lq\lq generalized Pisot 
numbers\rq\rq
as accumulation points of \lq\lq generalized Salem numbers\rq\rq.
Trees of degree $\geq 2$ lead to iterated accumulation
points.

\section{Counting Schr\"odinger operators 
according to values of the Jacobi symbol}\label{sectjacobi}

We describe a local invariant for computing the number of 
Schr\"odinger operators of a tree over a field
$\mathbb F_q$ of odd characteristic according to values of the Jacobi-symbol 
(corresponding to coefficients which are 
zero, non-zero squares or non-squares) at diagonal entries indexed by vertices.

Our formulae define in fact a coloured local invariant as
defined in Section \ref{sectcolmonoid} with five free parameters 
$\epsilon_v,q_v,s_v,x_v,y_v$ for each vertex $v$ 
and an additional global parameter $q$. The parameters $\epsilon_v,q_v$
(and the global parameter $q$)
are involved in the monoid-product which depends thus on the root-colour,
see Section \ref{sectcolmonoid}.

Some specializations of this invariant count invertible 
Schr\"odinger operators
with various restrictions (having for example only non-zero squares
on the diagonal, or zeros, non-zero squares, respectively non-squares
on selected subsets of diagonal entries).

Verifications are straightforward but tedious and are omitted.

We denote by $\mathcal O$ either the zero element of $\mathbb F_q$ or the element
$[0]$ of the group-algebra $\mathbb Q[\mathbb F_q]$ of the additive group 
$\mathbb F_q$. Similarly, $\mathcal S$ is either the set of 
all non-zero squares of the field $\mathbb F_q$ or the weighted sum 
$\frac{2}{q-1}\sum_{s\in \mathcal S}[s]$ of all non-zero squares in $\mathbb Q[\mathbb F_q]$, with equal weights summing up to $1$.
We define $\mathcal N$ analogously using non-squares.

Addition-rules for the elements $\mathcal S,\mathcal N$ (corresponding
to the product in the group-algebra $\mathbb Q[\mathbb F_q]$ of the 
additive group $(\mathbb F_q,+)$) 
are given by
$$\begin{array}{c|cc}
&\mathcal S&\mathcal N\\
\hline
\mathcal S&\frac{1+\epsilon}{q-1}\mathcal O+\frac{q-4-\epsilon}{2(q-1)}\mathcal S+\frac{q-\epsilon}{2(q-1)}\mathcal N&
\frac{1-\epsilon}{q-1}\mathcal O+\frac{q-2+\epsilon}{2(q-1)}(\mathcal S+\mathcal N)\\
\mathcal N&\frac{1-\epsilon}{q-1}\mathcal O+\frac{q-2+\epsilon}{2(q-1)}(\mathcal S+\mathcal N)&
\frac{1+\epsilon}{q-1}\mathcal O+\frac{q-\epsilon}{2(q-1)}\mathcal S+\frac{q-4-\epsilon}{2(q-1)}\mathcal N
\end{array}$$
where $\epsilon\in\pm 1$ satisfies $q\equiv \epsilon\pmod 4$.

Addition-rules in $\mathbb Q[\mathbb F_q]$ 
with $\mathcal O$ are of course given by $\mathcal O+X=X$ for 
$X\in \{\mathcal O,S\mathcal ,\mathcal N\}$. 

Multiplication-rules are given by 
$\mathcal O^2=\mathcal O,\mathcal O\cdot \mathcal S=\mathcal O\cdot 
\mathcal N=\mathcal O$,
$\mathcal S^2=\mathcal N^2=\mathcal S$ and $\mathcal S\mathcal N=\mathcal N$. In particular, the element $\mathcal S$ is a multiplicative 
unity.

Given a rooted tree $R$ and two subsets $\mathcal A,\mathcal B$ of a 
finite field $\mathbb F_q$, we use the conventions of Section 
\ref{sectmainexple} and
we denote by $\nu(\mathcal A,\mathcal B)$ the number of Schr\"odinger operators
with diagonals in $\mathbb F_q$ of the rooted tree $R$ having determinants 
of the form
$\alpha x+\beta$ with $\alpha\in \mathcal A$ and $\beta\in \mathcal B$.
We denote by $\mathcal O=\{0\}$ the zero element of $\mathbb F_q$ and 
by $\mathcal S$, respectively $\mathcal N$ the set of squares, respectively 
non-squares in $\mathbb F_q$ where we assume that $q$ is a power of an odd 
prime. We display all possible numbers $\nu(\mathcal A,\mathcal B)$
with $\mathcal A,\mathcal B\in\{\mathcal O,\mathcal S,\mathcal N\}$
of a rooted tree $R$ as a square-matrix with rows and columns 
indexed by $\mathcal O,\mathcal S$ and $\mathcal N$ by writing
$$\left(\begin{array}{ccc}
\nu(\mathcal O,\mathcal O)&\nu(\mathcal O,\mathcal S)&\nu(\mathcal O,\mathcal N)\\
\nu(\mathcal S,\mathcal O)&\nu(\mathcal S,\mathcal S)&\nu(\mathcal S,\mathcal N)\\
\nu(\mathcal N,\mathcal O)&\nu(\mathcal N,\mathcal S)&\nu(\mathcal N,\mathcal N)\end{array}\right).$$
The number $\nu(\mathcal O,\mathcal O)$ is in fact useless from our point 
of view but it allows for some numerical consistency checks.

The element
$$v=\left(\begin{array}{ccc}0&0&0\\
1&0&0\\
0&0&0\end{array}\right)$$
is the identity of monoid structures with products given by 
\begin{eqnarray*}
m_v(\left(\begin{array}{ccc}
a_1&b_1&c_1\\d_1&e_1&f_1\\g_1&h_1&i_1\end{array}\right),
\left(\begin{array}{ccc}
a_2&b_2&c_2\\d_2&e_2&f_2\\g_2&h_2&i_2\end{array}\right))
&=&\left(\begin{array}{ccc}A&B&C\\D&E&F\\G&H&I\end{array}\right)
\end{eqnarray*}
where
\begin{eqnarray*}
A&=&a_1a_2+a_1(b_2+c_2+d_2+e_2+f_2+g_2+h_2+i_2)\\
&&+(b_1+c_1+d_1+e_1+f_1+g_1+h_1+i_1)a_2+(b_1+c_1)(b_2+c_2)\\
B&=&b_1(d_2+e_2+f_2)+c_1(g_2+h_2+i_2)\\
&&+(d_1+e_1+f_1)b_2+(g_1+h_1+i_1)c_2\\
C&=&c_1(d_2+e_2+f_2)+b_1(g_2+h_2+i_2)\\
&&+(d_1+e_1+f_1)c_2+(g_1+h_1+i_1)b_2\\
D&=&d_1d_2+g_1g_2+\frac{1+\epsilon_v}{q-1}(e_1e_2+i_1i_2+f_1f_2+h_1h_2)\\
&&+\frac{1-\epsilon_v}{q-1}(e_1f_2+f_1e_2+h_1i_2+i_1h_2)\\
%+\frac{1+\epsilon_v}{q-1}(f_1f_2+h_1h_2)\\
E&=&d_1e_2+e_1d_2+g_1i_2+i_1g_2+\frac{q_v-4-\epsilon_v}{2(q-1)}(e_1e_2+i_1i_2)\\
&&+\frac{\tilde q_v-2+\epsilon_v}{2(q-1)}(e_1f_2+f_1e_2+h_1i_2+i_1h_2)+\frac{q_v-\epsilon_v}{2(q-1)}(f_1f_2+h_1h_2)\\
F&=&d_1f_2+f_1d_2+g_1h_2+h_1g_2+\frac{\tilde q_v-4-\epsilon_v}{2(q-1)}(f_1f_2+h_1h_2)\\
&&+\frac{q_v-2+\epsilon_v}{2(q-1)}(e_1f_2+f_1e_2+h_1i_2+i_1h_2)+\frac{\tilde q_v-\epsilon_v}{2(q-1)}(e_1e_2+i_1i_2)\\
G&=&d_1g_2+g_1d_2+\frac{1-\epsilon_v}{q-1}(e_1h_2+h_1e_2+f_1i_2+i_1f_2)\\
&&+\frac{1+\epsilon_v}{q-1}(e_1i_2+i_1e_2+f_1h_2+h_1f_2)\\
H&=&d_1h_2+h_1d_2+f_1g_2+g_1f_2+\frac{\tilde q_v-4-\epsilon_v}{2(q-1)}(f_1h_2+h_1f_2)\\
&&+\frac{q_v-2+\epsilon_v}{2(q-1)}(e_1h_2+h_1e_2+f_1i_2+i_1f_2)+\frac{\tilde q_v-\epsilon_v}{2(q-1)}(e_1i_2+i_1e_2)\\
I&=&d_1i_2+i_1d_2+e_1g_2+g_1e_2+\frac{q_v-4-\epsilon_v}{2(q-1)}(e_1i_2+i_1e_2)\\
&&+\frac{\tilde q_v-2+\epsilon_v}{2(q-1)}(e_1h_2+h_1e_2+f_1i_2+i_1f_2)+\frac{q_v-\epsilon_v}{2(q-1)}(f_1h_2+h_1f_2)\\
\end{eqnarray*}
and where $q_v,\tilde q_v$ satisfy
$$q_v+\tilde q_v=2q$$
and have thus only one degree of freedom at a given vertex since $q$ is
global (i.e. independent of vertices).

$\epsilon_v$ are also local variables.

It can be checked that $m_v$ defines for all $q,q_v,\epsilon_v$ (with $\tilde q_v=2q-q_v$) a commutative and associative product with identity $v$.

The reader should be aware that we use the same letter $v$ for a vertex, its colour and the 
identity of colour $v$ with respect to the monoid with product $m_v$.

Edge-operators do not depend on the final root-colour which is thus 
omitted. The index $v$ in our formulae denotes the root-vertex 
(or its colour) of
the rooted tree given as the argument. The edge-operator $e_v$ is given by
\begin{eqnarray*}
e_v(\left(\begin{array}{ccc}
a&b&c\\d&e&f\\g&h&i\end{array}\right))&=&
%&&=\left(\frac{q-1}{2}\right)^2
\left(\begin{array}{ccc}A&\frac{1+\epsilon_v}{2}B_++\frac{1-\epsilon_v}{2}B_-&
\frac{1+\epsilon_v}{2}C_++\frac{1-\epsilon_v}{2}C_-\\
D&\frac{1+\epsilon_v}{2}E_++\frac{1-\epsilon_v}{2}E_-&
\frac{1+\epsilon_v}{2}F_++\frac{1-\epsilon_v}{2}F_-\\
G&\frac{1+\epsilon_v}{2}H_++\frac{1-\epsilon_v}{2}H_-&
\frac{1+\epsilon_v}{2}I_++\frac{1-\epsilon_v}{2}I_-\end{array}\right)
\end{eqnarray*}
(using the conventions of Section \ref{sectmainexple}, 
except for a factor $((q-1)/2)^2$ 
corresponding to arbitrary non-zero square values on oriented edges,
when working with the specialization $q_v=q$ an odd prime power and $\epsilon_v\in\{\pm 1\}$ given by $q\equiv \epsilon_v\pmod 4$)
where 
\begin{eqnarray*}
A&=&\left(s_v+\frac{q_*-1}{2}(x_v+y_v)\right)a\\
B_+&=&s_vd+x_ve+y_vf\\
B_-&=&s_vg+x_vh+y_vi\\
C_+&=&s_vg+x_vi+y_vh\\
C_-&=&s_vd+x_vf+y_ve\\
\end{eqnarray*}
\begin{eqnarray*}
D&=&\left(s_v+\frac{q-1}{2}(x_v+y_v)\right)b\\
E_+&=&s_ve+\frac{q-1}{2}x_vd+\frac{q_v-4-\epsilon_v}{4}x_ve+\frac{\tilde q_v-2+\epsilon_v}{4}(x_vf+y_ve)+\frac{q_v-\epsilon_v}{4}y_vf\\
E_-&=&s_vh+\frac{q-1}{2}y_vg+\frac{\tilde q_v-4-\epsilon_v}{4}y_vh
+\frac{q_v-2+\epsilon_v}{4}(x_vh+y_vi)+\frac{\tilde q_v-\epsilon_v}{4}x_vi\\
F_+&=&s_vh+\frac{q-1}{2}y_vg+\frac{\tilde q_v-4-\epsilon_v}{4}y_vh
+\frac{q_v-2+\epsilon_v}{4}(x_vh+y_vi)+\frac{\tilde q_v-\epsilon_v}{4}x_vi\\
F_-&=&s_ve+\frac{q-1}{2}x_vd+\frac{q_v-4-\epsilon_v}{4}x_ve
+\frac{\tilde q_v-2+\epsilon_v}{4}(x_vf+y_ve)+\frac{q_v-\epsilon_v}{4}y_vf\\
G&=&\left(s_v+\frac{q-1}{2}(x_v+y_v)\right)c\\
H_+&=&s_vf+\frac{q-1}{2}y_vd+\frac{\tilde q_v-4-\epsilon_v}{4}y_vf
+\frac{q_v-2+\epsilon_v}{4}(x_vf+y_ve)+\frac{\tilde q_v-\epsilon_v}{4}x_ve\\
H_-&=&s_vi+\frac{q-1}{2}x_vg+\frac{q_v-4-\epsilon_v}{4}x_vi
+\frac{\tilde q_v-2+\epsilon_v}{4}(x_vh+y_vi)+\frac{q_v-\epsilon_v}{4}y_vh\\
I_+&=&s_vi+\frac{q-1}{2}x_vg+\frac{q_v-4-\epsilon_v}{4}x_vi
+\frac{\tilde q_v-2+\epsilon_v}{4}(x_vh+y_vi)+\frac{q_v-\epsilon_v}{4}y_vh\\
I_-&=&s_vf+\frac{q-1}{2}y_vd+\frac{\tilde q_v-4-\epsilon_v}{4}y_vf
+\frac{q_v-2+\epsilon_v}{4}(x_vf+y_ve)+\frac{\tilde q_v-\epsilon_v}{4}x_ve\\
\end{eqnarray*}
where we have $q_v+\tilde q_v=2q$, as above.

We consider moreover
\begin{eqnarray*}
f^s_v&=&D+\frac{1+\epsilon_v}{2}E_++\frac{1-\epsilon_v}{2}E_-+
\frac{1+\epsilon_v}{2}F_++\frac{1-\epsilon_v}{2}F_-\\
f^n_v&=&G+\frac{1+\epsilon_v}{2}H_++\frac{1-\epsilon_v}{2}H_-+
\frac{1+\epsilon_v}{2}I_++\frac{1-\epsilon_v}{2}I_-
\end{eqnarray*}
with $D,E_\pm,F_\pm,G,H_\pm,I_\pm$ as in the definition of $E_v$
and we set 
$$f_v=f_v^s+f_v^n\ .$$

We have then for arbitrary $q,q_v,\epsilon_v,s_v,x_v,y_v,q_w,\epsilon_w,s_w,x_w,y_w$ with $\tilde q_v=2q-q_v,\tilde q_w=2q-q_w$ the identity
\begin{eqnarray}\label{traceidjacobi}
f_v(m_v(A,e_w(B)))=f_w(m_w(e_v(A),B))
\end{eqnarray}
which holds for all $A,B$ and which corresponds 
to (\ref{formulecolmonoid}). We get thus a coloured local invariant. Observe that this identity is surprisingly general: The construction
of the invariant (enumeration of invertible Schr\"odinger operators
according to values of the Jacobi symbol on the diagonal) ensures it only for
$\epsilon_v=\epsilon$ with $\epsilon\in \{\pm 1\}$ such that $q
\equiv \epsilon\pmod 4$ for $q=q_v=\tilde q_v=q_w=\tilde q_w$ an odd prime-power.

Choosing an odd prime-power $q$ and setting $q_v=q$ and $\epsilon_v=\epsilon$
where $\epsilon\in\{\pm 1\}$ is given by $q\equiv \epsilon\pmod 4$,
the coefficient of a monomial
$$\prod_{v\in V(T)}u_v$$
with $u_v\in\{s_v,x_v,y_v\}$ of the final result
is a natural integer counting the number of invertible Schr\"odinger operators
with diagonal coefficient at a vertex $v$ equal $0$ if $u_v=s_v$,
respectively in the set of non-zero squares or non-squares if
$u_v=x_v$ or $u_v=y_v$.
Contributions from $f^s$ correspond to operators with determinant a non-zero
square and contributions from $f^n$ correspond to non-square determinants.

Replacing $f_v$ either by $f^s_v$ or by $f^n_v$ we lose in general 
identity (\ref{traceidjacobi}). It remains however valid for a few specializations. Example are:
\begin{itemize}
\item{} $\epsilon_v=\epsilon$ for $\epsilon\in\{\pm 1\}$,
$q_v=q$ (and $s_v,x_v,y_v$ arbitrary at each 
vertex). This case, with $\epsilon\equiv q\pmod 4$,
counts of course Schr\"odinger operators with determinant a non-zero 
square, respectively a non-square, of $\mathbb F_q$. 
\item{} $\epsilon_v=\epsilon,\ q_v=q,\ s_v=\frac{1-q}{2}(x_v+y_v)$.
\item{} $\epsilon_v=\epsilon,\ q_v=q,\ x_v=y_v=x$.
\item{} $\epsilon_v=\epsilon,\ q_v=q$ and $y_v=-x_v$ at each vertex $v$ 
of $T$.
\end{itemize}

\section{An identity in $\mathbb Z[\mathrm{SL}_2(\mathbb F_q)]$
}\label{sectgroup}
We compute the coefficients of
\begin{eqnarray}\label{nthpower}
\left(
\sum_{\mu\in M}\sum_{x\in\mathbb F_q}\left[\left(\begin{array}
{cc}x&\mu\\
-1/\mu&0\end{array}\right)\right]\right)^n\in\mathbb Z[\mathrm{SL}_2(\mathbb F_q)]
\end{eqnarray}
where $M$ is a subgroup of the multiplicative group of units
$\mathbb F_q^*$.
This allows to compute the sum over all coefficients corresponding to elements
of trace $2$ or of trace $-2$. For $M=\{1\}$ the trivial group, 
these sums encode the number of Schr\"odinger operators of determinant zero
over $\mathbb F_q$ for the $n$-cycle $C_n$.

\begin{rem}\label{remfactorization} 
The factorizations
\begin{eqnarray}\label{factidentity}
\sum_{\mu\in M}\sum_{x\in\mathbb F_q}\left[\left(\begin{array}
{cc}x&\mu\\
-1/\mu&0\end{array}\right)\right]&=&A(M)U=UA(M)
\end{eqnarray}
where 
$$A(M)=\sum_{\mu\in M}\left[\left(\begin{array}
{cc}\mu^{-1}&0\\0&\mu\end{array}\right)\right]\hbox{ and }
U=\sum_{x\in\mathbb F_q}\left[\left(\begin{array}
{cc}x&1\\
-1&0\end{array}\right)\right]$$
show that $M$ plays only a minor role in (\ref{nthpower}).

The element $A(M)$ satisfies $A(M)^2=mA(M)$ where $m$ is the number of elements 
of $M$. The element $\frac{1}{m}A(M)$ is thus
a non-central idempotent (of rank $\frac{q^3-q}{m}$) in $\mathbb Q[\mathrm{SL}_2(\mathbb F_q)]$. The spectrum (with multiplicities) of the linear 
endomorphism of $\mathbb Q[\mathrm{SL}_2]$ defined by $X\longmapsto UX$ 
can also easily be
recovered from our data. It consists of a subset of 
$\{-1,0,\pm \sqrt{q},\pm\sqrt{-q},q\}$, with $\pm\sqrt{-q}$ only occurring
if $-1\not\in M$.
\end{rem}

For $n\geq 1$ we consider the three sequences
\begin{eqnarray*}
\alpha_n&=&\left(q^{(n-1)\pmod 2}\right)\frac{q^{2\lfloor (n-1)/2\rfloor}-1}{q^2-1}
-\frac{q^{\lfloor (n-1)/2\rfloor}-1}{q-1}\\
\beta_n&=&\left(q^{(n-1)\pmod 2}\right)\frac{q^{2\lfloor (n-1)/2\rfloor}-1}{q^2-1}\\
\gamma_n&=&\left(q^{(n-1)\pmod 2}\right)\frac{q^{2\lfloor (n-1)/2\rfloor}-1}{q^2-1}
-\frac{q^{\lfloor (n-1)/2\rfloor}-1}{q-1}+\frac{q^{\lfloor (n-1)/2\rfloor}}{m}
\end{eqnarray*}
(with $\gamma_n$ depending on $m$, considered as a fixed constant)
with $(n-1)\pmod 2\in\{0,1\}$ equal $1$ if $n$ is even and zero otherwise
and with $\lfloor (n-1)/2\rfloor=\frac{n-2}{2}$ for even $n$
and $\lfloor (n-1)/2\rfloor=\frac{n-1}{2}$ for odd $n$.

For $q\geq 2$ a prime-power and for $m$ a natural integer dividing $q-1$ we have $\alpha_n\leq \beta_n<\gamma_n$ for all $n\geq 1$
with equality occurring only for $\alpha_1=\beta_1=0$ and $\alpha_2=\beta_2=0$.

Given a multiplicative subgroup $M$ of $m$ elements in the  
unit group $\mathbb F_q^*$, we set $\overline M=\mathbb F_q^*\setminus M$
if $-1\in M$, respectively $\overline M=\mathbb F_q^*
\setminus(+M\cup -M)$ if $-1\not \in M$.
Observe that $-1\not\in M$ if and only if $mq$ is odd. 
Given two subsets $B,D$ of $\mathbb F_q$, we denote by 
$\left(\begin{array}{c}B\\D\end{array}\right)$ the subset of all column-vectors
of $\mathbb F_q^2$ with first coordinate in $B$ and second coordinate in $D$.
We consider now the partition of all non-zero elements of $\mathbb F_q^2$
given by the four disjoint subsets 
$$
\left(\begin{array}{c}M\\0\end{array}\right),
\left(\begin{array}{c}\overline{M}\\0\end{array}\right),
\left(\begin{array}{c}\mathbb F_q\\M\end{array}\right),
\left(\begin{array}{c}\mathbb F_q\\\overline{M}\end{array}\right)$$
if $-1\in M$ (i.e. if $mq$ is even), respectively by the six disjoint subsets
$$\left(\begin{array}{c}+M\\0\end{array}\right),
\left(\begin{array}{c}-M\\0\end{array}\right),
\left(\begin{array}{c}\overline{M}\\0\end{array}\right),
\left(\begin{array}{c}\mathbb F_q\\+M\end{array}\right),
\left(\begin{array}{c}\mathbb F_q\\-M\end{array}\right),
\left(\begin{array}{c}\mathbb F_q\\\overline{M}\end{array}\right)
$$
if $-1\not\in M$ (i.e. if $mq$ is odd) where $0$ denotes of course 
the singleton $\{0\}\subset \mathbb F_q$.

For a fixed multiplicative 
subgroup $M$ of $m$ elements in $\mathbb F_q^*$ and $n\geq 1$, we consider the four, respectively six, 
rational sequences, named by the parts of the above partition of 
non-zero column-vectors in $\mathbb F_q^2$, given by table (\ref{formulaseqvalues}).
\begin{eqnarray}\label{formulaseqvalues}
\begin{array}{||c||c|c|c|c||}
\hline\hline
n\pmod 4:&0&1&2&3\\
\hline\hline
\left(\begin{array}{c}M\\0\end{array}\right)_n&
\beta_n&\gamma_n&\beta_n&\gamma_n\\
\hline
\left(\begin{array}{c}+M\\0\end{array}\right)_n&
\beta_n&\gamma_n&\beta_n&\alpha_n\\
\hline
\left(\begin{array}{c}-M\\0\end{array}\right)_n&
\beta_n&\alpha_n&\beta_n&\gamma_n\\
\hline
\left(\begin{array}{c}\overline{M}\\0\end{array}\right)_n&
\beta_n&\alpha_n&\beta_n&\alpha_n\\
\hline
\left(\begin{array}{c}\mathbb F_q\\M\end{array}\right)_n&
\gamma_n&\beta_n&\gamma_n&\beta_n\\
\hline
\left(\begin{array}{c}\mathbb F_q\\+M\end{array}\right)_n&
\gamma_n&\beta_n&\alpha_n&\beta_n\\
\hline
\left(\begin{array}{c}\mathbb F_q\\-M\end{array}\right)_n&
\alpha_n&\beta_n&\gamma_n&\beta_n\\
\hline
\left(\begin{array}{c}\mathbb F_q\\\overline{M}\end{array}\right)_n&
\alpha_n&\beta_n&\alpha_n&\beta_n\\
\hline\hline
\end{array}
\end{eqnarray}
Values for $\left(\begin{array}{c}M\\0\end{array}\right)_n$
and $\left(\begin{array}{c}\mathbb F_q\\M\end{array}\right)_n$ 
apply only if $-1\in M$ and depend only on the parity of $n$.
Values for $\left(\begin{array}{c}\overline{M}\\0\end{array}\right)_n$ 
and $\left(\begin{array}{c}\mathbb F_q\\\overline{M}\end{array}\right)_n$
apply whether or not $-1$ is in $M$ and depend also only on the parity of $n$.
The remaining values involving $+M$ or $-M$ apply only if $-1\not\in M$
and depend on $n$ modulo $4$.

For a non-zero vector $\left(\begin{array}{c}b\\d\end{array}\right)$
and for $n\geq1$
(and for a given fixed subgroup $M\subset \mathbb F_q^*$ of $m$ elements) we 
set $\left(\begin{array}{c}b\\d\end{array}\right)_n=\left(\begin{array}{c}B\\D\end{array}\right)_n$
if $b\in B$ and $d\in D$ with $B,D\in \{M,\overline M,0,\mathbb F_q\}$
in the case where $-1\in M$,
respectively with $B,D\in \{+M,-M,\overline M,0,\mathbb F_q\}$
in the case where $-1\not\in M$. Observe that 
$\left(\begin{array}{c}b\\d\end{array}\right)_n\in\{\alpha_n,\beta_n,\gamma_n\}$
for all $n\geq 1$. We have now the following result:

\begin{thm}\label{thmgroup} For all integers $n\geq 1$ we have the identity
\begin{eqnarray*}
%\label{groupidentity}
\left(\sum_{\mu\in M}\sum_{x\in\mathbb F_q}\left[\left(\begin{array}{cc}x&\mu\\-1/\mu&0\end{array}
\right)\right]\right)^n&=&
m^n\sum_{\mathrm{SL}_2(\mathbb F_q)}
\left(\begin{array}{c}b\\d\end{array}\right)_n\
\left[\left(\begin{array}{cc}a&b\\c&d\end{array}\right)\right]
\end{eqnarray*}
in the group ring $\mathbb Z[\mathrm{SL}_2(\mathbb F_q)]$
where $m$ is the number of elements in a subgroup $M$ of $\mathbb F_q^*$
and where $\left(\begin{array}{c}b\\d\end{array}\right)_n$ are
as above
(the last sum is of course 
over all $q^3-q$ elements 
$\left(\begin{array}{cc}a&b\\c&d\end{array}\right)$ of 
$\mathrm{SL}_2(\mathbb F_q)$).
\end{thm}

\begin{rem} Theorem \ref{thmgroup} can easily be adapted to
$\mathrm{PSL}_2(\mathbb F_q)$.
\end{rem}

The following Lemma is the main ingredient for proving Theorem \ref{thmgroup}:
\begin{lem}\label{lemrecidabc} 
The three sequences $\alpha_n,\beta_n,\gamma_n,\ n\geq 1$ satisfy 
the identities
\begin{eqnarray*}
\alpha_{2n}&=&\alpha_{2n-1}+(q-1)\beta_{2n-1},\\
\alpha_{2n+1}&=&q\alpha_{2n},\\
\beta_n&=&q\beta_{n-1}+\left(n\pmod 2\right),\\
&=&(q-1-m)\alpha_{n-1}+\beta_{n-1}+m\gamma_{n-1}-\left((n-1)\pmod 2\right),\\
\gamma_{2n}&=&(q-1)\beta_{2n-1}+\gamma_{2n-1},\\
\gamma_{2n+1}&=&q\gamma_{2n}.\\
\end{eqnarray*}
\end{lem}
We leave the straightforward but tedious verifications to the
reader.\hfill$\Box$

\begin{rem} The two recursive identities for $\beta_n$ yield
\begin{eqnarray*}(q-1-m)\alpha_n+m\gamma_n&=&(q-1)\beta_n+1\\
%(q-2)\alpha_{2n}+\gamma_{2n-1}&=&(q-1)^2\beta_{2n-1}+1.
\end{eqnarray*}
which shows linear dependency of the constant sequence $1,1,1,\dots$ from 
the three sequences $\alpha_n,\beta_n,\gamma_n$.
\end{rem}

\noindent{\bf Proof of Theorem \ref{thmgroup}} The result holds 
by (\ref{formulaseqvalues}) for $n=1$ 
since $\alpha_1=\beta_1=0$ and $\gamma_1=\frac{1}{m}$.

For $n\geq 2$, we have the recursive formula
\begin{eqnarray*}
m^n\left(\begin{array}{c}b\\d\end{array}\right)_n&=&m^{n-1}\sum_{\mu\in M}\sum_{x\in
\mathbb F_q}\left(
\begin{array}{cc}0&-\mu\\1/\mu&x\end{array}\right)\left(\begin{array}{c}
b\\d\end{array}\right)_{n-1}
\end{eqnarray*}
and induction implies
\begin{eqnarray}\label{recformcoeffn}
m\left(\begin{array}{c}b\\d\end{array}\right)_n&=&\sum_{\mu\in M}\sum_{x\in
\mathbb F_q}\left(\begin{array}{c}
-\mu d\\b/\mu+xd\end{array}\right)_{n-1}.
\end{eqnarray}
Elementary properties of finite fields imply now
the recursive formulae
\begin{eqnarray*}
%\label{recformeven}
\left(\begin{array}{c}M\\0\end{array}\right)_n
&=&q\left(\begin{array}{c}\mathbb F_q\\M\end{array}\right)_{n-1}\\
\left(\begin{array}{c}\overline{M}\\0\end{array}\right)_n
&=&q\left(\begin{array}{c}\mathbb F_q\\\overline{M}\end{array}\right)_{n-1}\\
\left(\begin{array}{c}\mathbb F_q\\M\end{array}\right)_n
&=&\left(\begin{array}{c}M\\0\end{array}\right)_{n-1}+m\left(\begin{array}{c}\mathbb F_q\\M\end{array}\right)_{n-1}+(q-1-m)\left(\begin{array}{c}\mathbb F_q\\\overline{M}\end{array}\right)_{n-1}\\
\left(\begin{array}{c}\mathbb F_q\\\overline{M}\end{array}\right)_n
&=&\left(\begin{array}{c}\overline{M}\\0\end{array}\right)_{n-1}+m\left(\begin{array}{c}\mathbb F_q\\M\end{array}\right)_{n-1}+(q-1-m)\left(\begin{array}{c}\mathbb F_q\\\overline{M}\end{array}\right)_{n-1}
\end{eqnarray*}
if $-1\in M$ and the recursive identities
\begin{eqnarray*}
\left(\begin{array}{c}\pm M\\0\end{array}\right)_n
&=&q\left(\begin{array}{c}\mathbb F_q\\\pm M\end{array}\right)_{n-1}\\
\left(\begin{array}{c}\overline{M}\\0\end{array}\right)_n
&=&q\left(\begin{array}{c}\mathbb F_q\\\overline{M}\end{array}\right)_{n-1}\\
\left(\begin{array}{c}\mathbb F_q\\\pm M\end{array}\right)_n
&=&\left(\begin{array}{c}\mp M\\0\end{array}\right)_{n-1}+T_{n-1}\\
\left(\begin{array}{c}\mathbb F_q\\\overline{M}\end{array}\right)_n
&=&\left(\begin{array}{c}\overline{M}\\0\end{array}\right)_{n-1}+T_{n-1}
\end{eqnarray*}
where
$$
T_{n-1}=m\left(\begin{array}{c}\mathbb F_q\\+M\end{array}\right)_{n-1}+
m\left(\begin{array}{c}\mathbb F_q\\-M\end{array}\right)_{n-1}+(q-1-2m)
\left(\begin{array}{c}\mathbb F_q\\ \overline{M}\end{array}\right)_{n-1}
$$
if $-1\not\in M$.
Replacing all expressions by their values given by 
(\ref{formulaseqvalues}) we check that all these
expressions boil down to equalities
of Lemma \ref{lemrecidabc}.\hfill$\Box$

\begin{rem} Spectral calculus gives a different proof of 
Theorem \ref{thmgroup}.
\end{rem}

\subsection{Traces}

We denote by $m^nS_\tau(n)$ the sum of all coefficients
in (\ref{nthpower}) 
corresponding to elements of trace $\tau\in\mathbb F_q$. 
We are only interested in the values of $S_2(n)$ and 
$S_{-2}(n)$. Observe that $\mathrm{SL}_2(\mathbb F_q)$ contains 
exactly $q^2$ elements of trace $2$: The identity and 
all (other) $q^2-1$ unipotent elements. Multiplication by $-1$ 
induces of course a bijection between elements of 
trace $\tau$ and elements of trace $-\tau$.

We have the following result:

\begin{prop}\label{propsumtraces}
If $n$ is even we have
\begin{align}
%\label{formsumtraceeven}
S_{\pm 2}(n)=&\ q^3\frac{q^{n-2}-1}{q^2-1}-(q^2-q+1)\frac{q^{(n-2)/2}-1}{q-1}\\
&+\frac{m(q-1)+q\kappa(\pm 2,n)}{m}q^{(n-2)/2}\nonumber
\end{align}
with $\kappa(\pm 2,n)=1$ if $-1\in M$ and with $\kappa(\pm 2,n)$ given 
by
$$\begin{array}{||c||c|c||}
\hline\hline
&n\equiv 0\pmod 4&n\equiv 2\pmod 4\\
\hline\hline
\kappa(2,n)&1&0\\
\hline
\kappa(-2,n)&0&1\\
\hline\hline
\end{array}$$
if $-1\not\in M$.
 
If $n$ is odd, we have 
\begin{eqnarray}\label{formsumtracesodd}
S_2(n)=S_{-2}(n)&=&\frac{q^{n+1}-1}{q^2-1},
\end{eqnarray}
independently of $M$.
\end{prop}

\noindent{\bf Proof}
Counting the number of matrices of trace $2$, respectively $-2$, 
in every possible class we get
\begin{eqnarray*}
S_{\pm2}(n)&=&m\left(\begin{array}{c}M\\0\end{array}\right)_n+(q-1-m)
\left(\begin{array}{c}\overline M\\0\end{array}\right)_n\\
&&+(m(q-1)+q)\left(\begin{array}{c}\mathbb F_q\\M\end{array}\right)_n+(q-1-m)(q-1)
\left(\begin{array}{c}\mathbb F_q\\\overline M\end{array}\right)_n.
\end{eqnarray*}
if $-1\in M$ and 
\begin{eqnarray*}
S_{\pm 2}(n)&=&m\left(\begin{array}{c}+M\\0\end{array}\right)_n+m\left(\begin{array}{c}-M\\0\end{array}\right)_n+(q-1-2m)
\left(\begin{array}{c}\overline M\\0\end{array}\right)_n\\
&&+(m(q-1)+q)\left(\begin{array}{c}\mathbb F_q\\\pm M\end{array}\right)_n+m(q-1)\left(\begin{array}{c}\mathbb F_q\\\mp M\end{array}\right)_n\\
&&+(q-1-2m)(q-1)
\left(\begin{array}{c}\mathbb F_q\\\overline M\end{array}\right)_n.
\end{eqnarray*}
if $-1\not\in M$.

If $-1\in M$, the above expression for $S_{\pm 2}(n)$ amounts to 
\begin{eqnarray*}
&&(q-1)\beta_n+(m(q-1)+q)\gamma_n+(q-1-m)(q-1)\alpha_n\\
&=&q^3\frac{q^{n-2}-1}{q^2-1}-(q^2-q+1)\frac{q^{(n-2)/2}-1}{q-1}+
(m(q-1)+q)\frac{q^{(n-2)/2}}{m}
\end{eqnarray*}
if $n$ is even and to 
\begin{eqnarray*}
m\gamma_n+(q-1-m)\alpha_n+(q^2-q+1)\beta_n&=&\frac{q^{n+1}-1}{q^2-1}\\
%&=&q^2\frac{q^{n-1}-1}{q^2-1}-(q-1)\frac{q^{(n-1)/2}-1}{q-1}+q^{(n-1)/2}\\
%&=&q^2\frac{q^{n-1}-1}{q^2-1}+1
\end{eqnarray*}
if $n$ is odd.

If $-1\not\in M$ and $n$ even the value of $S_\tau(n)$ with $\tau\in\{\pm 2\}$
equals 
\begin{eqnarray*}
&(q-1)\beta_n+(q-1-m)(q-1)\alpha_n+m(q-1)\gamma_n+q\eta_n(\tau)
\end{eqnarray*}
with $\eta_n(\tau)$ given by
$$\begin{array}{||c||c|c||}
\hline\hline
&n\equiv 0\pmod 4&n\equiv 2\pmod 4\\
\hline\hline
\eta_n(2)=&\gamma_n&\alpha_n\\
\hline
\eta_n(-2)=&\alpha_n&\gamma_n.\\
\hline\hline
\end{array}$$
For $\eta_n(\tau)=\alpha_n$ we have
$$q^3\frac{q^{n-2}-1}{q^2-1}-(q^2-q+1)\frac{q^{(n-2)/2}-1}{q-1}+(q-1)q^{(n-2)/2}$$
and for $\eta_n(\tau)=\gamma_n$ we get
$$q^3\frac{q^{n-2}-1}{q^2-1}-(q^2-q+1)\frac{q^{(n-2)/2}-1}{q-1}
+(m(q-1)+q)\frac{q^{(n-2)/2}}{m}.$$

For $-1\not\in M$ and $n$ odd, the common value $S_2(n)=S_{-2}(n)$ is given by 
\begin{eqnarray*}
m(\alpha_n+\gamma_n)+(q-1-2m)\alpha_n+(q^2-q+1)\beta_n
&=&\frac{q^{n+1}-1}{q^2-1}.
\end{eqnarray*}
This ends the proof.\hfill$\Box$

\section{Proof of Theorem \ref{thmcycles}}\label{sectproofthmcycles}

We denote by $I_n$ the graph having vertices $1,\dots,n$. Consecutive integers
represent adjacent vertices.
For $q$ a fixed prime-power and $n\geq 1$, we write 
$\nu\left(\begin{array}{cc}a&b\\c&d\end{array}\right)_n$ 
for the number of Schr\"odinger operators $M$ over $\mathbb F_q$ for $I_n$ 
such that 
$$
\begin{array}{ll}a=\det(M)&b=\det(M(n;n))\\ c=-\det(M(1;1))&d=-\det(M(1,n;1,n))
\end{array}$$
where $M(i;j)$ respectively $M(i_1,i_2;j_1,j_2)$ denotes the submatrix 
of $M$ obtained by deleting line(s) $i_*$ 
and row(s) $j_*$.

Initial values for $\nu\left(\begin{array}{cc}a&b\\c&d\end{array}\right)_1$
are given by 
$\nu\left(\begin{array}{cc}x&1\\-1&0\end{array}\right)_1=1$ and 
$\nu\left(\begin{array}{cc}a&b\\c&d\end{array}\right)_1=0$ if 
$(b,c,d)\not=(1,-1,0)$.
Expanding the determinant of a Schr\"odinger operator for $I_{n+1}$
with a first diagonal coefficient $x$ along the first row shows the
recursion
$$\nu\left(\begin{array}{cc}a&b\\c&d\end{array}\right)_{n+1}=
\sum_{x\in\mathbb F_q}\nu\left(\begin{array}{cc}-c&-d\\a+xc&b+xd\end{array}\right)_n.$$
The matrix identity
$$\left(\begin{array}{cc}a&b\\c&d\end{array}\right)=
\left(\begin{array}{cc}x&1\\-1&0\end{array}\right)
\left(\begin{array}{cc}-c&-d\\a+xc&b+xd\end{array}\right)$$
implies 
now the identity
$$\left(\sum_{x\in\mathbb F_q}\left[\left(\begin{array}{cc}x&1\\-1&0\end{array}\right)\right]\right)^n=\sum_{\mathrm{SL}_2(\mathbb F_q)}
\nu\left(\begin{array}{cc}a&b\\c&d\end{array}\right)_n\left[\left(\begin{array}{cc}a&b\\c&d\end{array}\right)\right]$$
in the group-ring $\mathbb Z[\mathrm{SL}_2(\mathbb F_q)]$.
This correspond of course to the case $M=1$
in (\ref{nthpower}) or in Theorem \ref{thmgroup}, with $-1\not\in M$ except if 
$\mathbb F_q$ is of characteristic $2$.

Given a Schr\"odinger operator $M$ for $I_n$ we denote by $\tilde M$
the Schr\"odinger operator of the $n$-cycle (obtained 
by joining the first and last vertex of $I_n$ with an additional edge)
with the same diagonal coefficients.
Thus, $\tilde M$ is obtained from $M$ by adding two non-zero coefficients $1$ at the upper-right and lower-left corner of $M$.
Denoting by $x_1$ the first diagonal coefficient of $M$ or $\tilde M$
and expanding the determinant of $\tilde M$ along the first row and perhaps 
subsequently along the first column, 
we get 
\begin{eqnarray*}
\det(\tilde M)&=&
x_1\det(M(1;1))-\det(M(1,2;1,2))\\
&&-(-1)^n\Big(\det(M(1,n;1,2))+\det(M(1,2;1,n))\Big)\\
&&-\det(M(1,n;1,n)).
\end{eqnarray*}
The identities 
\begin{eqnarray*}&x_1\det(M(1;1))-\det(M(1,2;1,2))=\det(M)=a,&\\
&\det(M(1,n;1,2))=\det(M(1,2;1,n))=1,&\\
&-\det(M(1,n;1,n))=c&
\end{eqnarray*}
with $\left(\begin{array}{cc}a&b\\c&d\end{array}\right)$
the matrix associated to $M$ as above yield
$$\det(\tilde M)=a+c-2(-1)^n.$$
The number of non-invertible Schr\"odinger operators for $C_n$
is thus the total sum $S_{2(-1)^n}(n)$ of coefficients 
associated to matrices of trace $2(-1)^n$
in $$\left(\sum_{x\in\mathbb F_q}\left[\left(\begin{array}{cc}x&1\\-1&0\end{array}\right)\right]\right)^n.$$
This shows that $q^n-S_{2(-1)^n}(n)$ is the number of invertible Schr\"odinger operators for $C_n$. Proposition \ref{propsumtraces}
with $m=1$ and $M=\{1\}$ the trivial subgroup of $\mathbb F_q$
gives the values for $S_{2(-1)^n}$. Observe that $-1\in M$ 
if $q$ is  a power of $2$ and $-1\not\in M$
if $q$ is an odd prime-power.
The easy identities
\begin{eqnarray*}
&&q^{2n}-\frac{(1-q^{2n})(1-q^{2n+1})}{1-q^2}\\
&=&q^3\frac{q^{4n-2}-1}{q^2-1}-(q^2-q+1)\frac{q^{2n-1}-1}{q-1}+(2q-1)q^{2n-1}
\end{eqnarray*}
and 
\begin{eqnarray*}
&&-\frac{(1-q^{2n+1})(1-q^{2n+2})}{1-q^2}\\
&=&q^3\frac{q^{4n}-1}{q^2-1}-(q^2-q+1)\frac{q^{2n}-1}{q-1}+(q-1)q^{2n}
\end{eqnarray*}
end the proof. \hfill$\Box$

\begin{rem} Since $S_2(n)=S_{-2}(n)$ for odd $n$, the number of invertible
Schr\"odinger operators for $C_n$ is also given by $q^n-S_2(n)$.
\end{rem}

\section{Proof for complete graphs}\label{sectproofcompletegraphs}

\subsection{Simple stars} 
The number of invertible Schr\"odinger operators of simple stars
(finite graphs with at most one non-leaf) is a crucial ingredient 
for proving Theorem \ref{thmcpltegraph}.

We denote in this subsection by $R_n$ the rooted graph consisting of 
a central root adjacent to $n-1$ leaves. Writing $S_r(R_n)=(a_n,b_n)$
with $S_r(R_n)$ defined as in Section \ref{sectmainexple},
we have
%$$a_n=(q-1)^{n-1}\hbox{ and }b_n=(n-1)(q-1)^{n-2}$$
$$S_{R_n}=\left((q-1)^{n-1},(n-1)(q-1)^{n-2}\right)$$
as can be checked using 
\begin{eqnarray*}
(a_{n+1},b_{n+1})&=&m((a_n,b_n),e(v))=m((a_n,b_n),(q-1,1))\\
&=&((q-1)a_n,a_n+(q-1)b_n)\ .
\end{eqnarray*}
We have thus
\begin{eqnarray}\label{Sstarn}
S_{*_n}&=&f(a_n,b_n)=\left(q^2+(n-3)q+1\right)(q-1)^{n-2}
\end{eqnarray}
for the non-rooted star $*_n$ underlying $R_n$ given by a central vertex
of degree $n-1$ surrounded by $n-1$ leaves.

\subsection{Proof of Theorem \ref{thmcpltegraph}}

We consider an invertible Schr\"odinger operator $M$ for 
the star $*_n$ with a central vertex of degree $n-1$ surrounded by $n-1$ leaves.
If the diagonal coefficient $\lambda$ of the central vertex is different
from $-1$, we get an invertible Schr\"odinger operator of the complete graph 
$K_n$ on $n$ vertices
by adding the first row, corresponding to the central vertex of $*_n$, 
of $M$ to all other rows and by dividing the first
column of the resulting matrix by $1+\lambda$.  This construction can be reversed, as easily seen
on the following illustration:
\begin{eqnarray*}
&\left(\begin{array}{cccccc}
\lambda&1&1&\dots&1\\
1&a_2&0&\dots&0\\
1&0&a_3&\dots&0\\
\vdots&&&\ddots&\vdots\\
1&0&0&\dots&0\\
1&0&0&\dots&a_n\end{array}\right)
\leftrightarrow
\left(\begin{array}{ccccc}
\lambda&1&1&\dots&1\\
1+\lambda&1+a_2&1&\dots&1\\
1+\lambda&1&1+a_3&\dots&1\\
\vdots&&&\ddots&\vdots\\
1+\lambda&1&1&\dots&1+a_n\\
\end{array}\right)&\\
&\leftrightarrow
\left(\begin{array}{cccccc}
\mu&1&1&\dots&1&1\\
1&1+a_2&1&\dots&1&1\\
1&1&1+a_3&\dots&1&1\\
\vdots&&&\ddots&&\vdots\\
1&1&1&\dots&a_{n-1}&1\\
1&1&1&\dots&1&1+a_n\end{array}\right).&
\end{eqnarray*}
Every invertible Schr\"odinger operator of $K_n$ with first diagonal
coefficient $\mu=\lambda/(1+\lambda)$ different from $1$ is of this form
by taking $\lambda=\frac{\mu}{1-\mu}$.

Schr\"odinger operators of $K_n$ with first coefficient $1$ are invertible 
if and only if they have a first diagonal coefficient 
equal to $1$ and all $n-1$ remaining diagonal coefficients different from $1$.
There are thus $(q-1)^{n-1}$ such matrices.

The number $\kappa_n$ of invertible Schr\"odinger operators for the complete
graph $K_n$ (over a fixed 
finite field $\mathbb F_q$) is thus given by
\begin{eqnarray}\label{formulakappan}
\kappa_n&=&S_{*_n}(q)+(q-1)^{n-1}-s_n
\end{eqnarray}
where $S_{*_n}$ is given by (\ref{Sstarn}) and 
where $s_n$ denotes the number of invertible Schr\"odinger operators of 
$*_n$ with $-1$ as the diagonal entry 
corresponding to the central vertex of the simple star $*_n$
consisting of a central vertex of degree $n-1$ adjacent to $n-1$ leaves.

The number $s_n$ can be computed as follows:
First observe that an invertible Schr\"odinger operator of $*_n$ 
has at most a unique diagonal coefficient which is zero. 
The contribution of such matrices to $s_n$, given by
\begin{eqnarray}\label{numberstarswithzero}
&(n-1)(q-1)^{n-2},&
\end{eqnarray}
is easy to establish.

Matrices contributing to $s_n$ having only non-zero diagonal entries 
are in bijection with solutions $(b_2,\dots,b_n)\in(\mathbb F_q^*)^{n-1}$
of 
$$-b_2\cdots b_n(1+\frac{1}{b_2}+\dots+\frac{1}{b_n})\not=0\ .$$
We denote by $\tilde s_n$ the number of 
such solutions $(b_2,\dots,b_n)\in (\mathbb F_q^*)^{n-1}$.

Choosing an arbitrary non-zero element $x_1\in\mathbb F_q^*$ and setting 
$x_2=\frac{x_1}{b_2},\dots,x_n=\frac{x_1}{b_n}$ we have $(q-1)\tilde s_n=\beta_n$
where $\beta_n$ counts the number of solutions $(x_1,\dots,x_n)\in(\mathbb F_q^*)^n$ of the inequality $\sum_{i=1}^n x_i\not=0$.
We denote similarily by $\alpha_n$ the number of solutions $(x_1,\dots,x_n)\in(\mathbb F_q^*)^n$ of the equality $\sum_{i=1}^n x_i=0$.

We have $\alpha_0=1,\ \beta_0=0$ and the recursive formulae
\begin{eqnarray*}
\alpha_n&=&\beta_{n-1},\\
\beta_n&=&(q-1)\alpha_{n-1}+(q-2)\beta_{n-1}.
\end{eqnarray*}

\begin{lem}\label{lembetan} We have
$$\beta_n=\frac{(q-1)^{n+1}+(-1)^{n+1}(q-1)}{q}.$$
\end{lem}

\noindent{\bf Proof} The sequence $\beta_n$ has the recursive definition
$\beta_0=0,\beta_1=q-1$ and $\beta_n=(q-2)\beta_{n-1}+(q-1)\beta_{n-2}$
for $n\geq 2$. An induction on $n$ based on the trivial identity
\begin{eqnarray*}
&&(q-1)^{n+1}+(-1)^{n+1}(q-1)\\
&=&(q-2)\left((q-1)^{n}+(-1)^{n}(q-1)\right)\\
&&+(q-1)\left((q-1)^{n-1}+(-1)^{n-1}(q-1)\right)
\end{eqnarray*}
ends the proof.\hfill$\Box$
\begin{rem} Lemma \ref{lembetan} follows also easily from the fact that 
the characteristic polynomial $X^2-(q-2)X-(q-1)$ 
of the recursion defining $\beta_n$ has roots $q-1$ and $-1$.
\end{rem}
 (\ref{numberstarswithzero}) and the equality $(q-1)\tilde s_n=\beta_n$
already mentionned
imply $s_n=(n-1)(q-1)^{n-2}+\frac{1}{q-1}\beta_n$. Lemma \ref{lembetan}
shows thus that
\begin{eqnarray}\label{formulastarn}
s_n&=&(n-1)(q-1)^{n-2}+\frac{(q-1)^{n}-(-1)^n}{q}.
\end{eqnarray}
Identities (\ref{Sstarn}), (\ref{formulakappan}) and (\ref{formulastarn}),
taken together, show that $K_n$ has
\begin{eqnarray*}
&&(q^2+(n-3)q+1)(q-1)^{n-2}+(q-1)^{n-1}\\
&&-\left((n-1)(q-1)^{n-2}+\frac{(q-1)^{n}-(-1)^n}{q}\right)\\
%&=&(q-1)^{n-2}\left(q^2+(n-3)q+1+q-1-(n-1)\right)+\frac{(-1)^n-(q-1)^n}{q}\\
&=&n(q-1)^{n-1}+\frac{(q-1)^{n-2}\left(q^3-2q^2+q-(q-1)^2\right)+(-1)^n}{q}\\
&=&n(q-1)^{n-1}+\frac{(q-1)^{n+1}+(-1)^n}{q}
\end{eqnarray*}
invertible Schr\"odinger operators over $\mathbb F_q$.\hfill$\Box$

\section{Final remarks}\label{sectcomplements}

\subsection{Generalizations}

It is of course possible to define Schr\"odinger operators for 
arbitrary (perhaps oriented) simple graphs and to count invertible 
Schr\"odinger operators over finite fields. I ignore if there
is an efficient way for computing the corresponding numbers.

Enumerating for example all Schr\"odinger operators over 
$\mathbb F_3$ and $\mathbb F_5$ for the Petersen graph
(obtained by identifying opposite points of the $1-$skeleton
of the dodecahedron) we get $q^{10}-q^9+q^8$ invertible
Schr\"odinger operators in both cases. This formula fails however
for $\mathbb F_2$ for which no invertible Schr\"odinger
operators exists.

A second notion, closely related to Schr\"odinger operators and 
used for example in Proposition \ref{propinvTmatrix},
is to look at the set of all 
invertible matrices with off-diagonal support defining a given
graph (diagonal elements are arbitrary). In the case of 
unoriented graphs, one can moreover require matrices to be symmetric.

For trees, both definitions are essentially identical (up to a factor 
$(q-1)^*$).

\subsection{Counting points over finite fields for algebraic varieties
over $\overline{\mathbb Q}$}

Our main problem, counting invertible Schr\"odinger operators, is 
of course a particular case of counting points over finite fields
on algebraic varieties defined over $\mathbb Z$ (or more generally
over $\overline{\mathbb Q}$). Such problems are in general difficult,
see e.g. the monograph \cite{Se} devoted to such questions.

The main problem is of course the question if the behaviour of 
invertible (or equivalently, non-invertible) Schr\"odinger operators
of graphs is simpler. 
All our examples give rise to polynomials (depending perhaps
on the parity of the characteristic) enumerating invertible 
Schr\"odinger operators over finite fields. Does this fail for
some finite graph $G$ or do there always exist polynomials 
depending on $q\pmod{N_G}$ for some natural 
integer $N_G$ evaluating to the number of invertible Schr\"odinger operators
over $\mathbb F_q$ for $G$?

\subsection{Acknowledgements}

I thank E. Peyre and P. De la Harpe for interesting discussions and remarks.

\noindent Roland BACHER, Univ. Grenoble Alpes, Institut 
Fourier (CNRS UMR 5582), 38000 Grenoble, France.

\noindent e-mail: Roland.Bacher@ujf-grenoble.fr

\end{document}